\documentclass[12pt]{article}%
\usepackage{amsmath}
\usepackage{amsfonts}
\usepackage{amssymb}
\usepackage{graphicx}%
\providecommand{\U}[1]{\protect\rule{.1in}{.1in}}
\providecommand{\U}[1]{\protect\rule{.1in}{.1in}}
\newtheorem{theorem}{Theorem}

\begin{document}

\title{Choreographies in the $n$-vortex problem}
\author{Renato Calleja \thanks{IIMAS, Universidad Nacional Aut\'{o}noma de M\'{e}xico,
calleja@mym.iimas.unam.mx}, Eusebius Doedel \thanks{Computer Science,
Concordia University, Montreal, Canada, doedel@cs.concordia.ca}, Carlos
Garc\'{\i}a-Azpeitia \thanks{Facultad de Ciencias, Universidad Nacional
Aut\'{o}noma de M\'{e}xico, cgazpe@ciencias.unam.mx}}
\maketitle

\begin{abstract}
We consider the equations of motion of $n$ vortices of equal circulation in
the plane, in a disk and on a sphere. The vortices form a polygonal
equilibrium in a rotating frame of reference. We use numerical continuation in
a boundary value setting to determine the Lyapunov families of periodic orbits
that arise from the polygonal relative equilibrium. When the frequency of a
Lyapunov orbit and the frequency of the rotating frame have a rational
relationship then the orbit is also periodic in the inertial frame. A dense
set of Lyapunov orbits, with frequencies satisfying a diophantine equation,
corresponds to choreographies of the $n$ vortices. We include numerical
results for all cases, for various values of $n$, and we provide key details
on the computational approach.

\end{abstract}

\section{Introduction}

\label{sec:intro}
A vortex point in an inviscid and incompressible two-di\-men\-sio\-nal fluid
corresponds to a singular vorticity concentrated at a point with constant
circulation, by using the Dirac delta function for the vorticity in the Euler
equation. This vorticity induces a singular velocity field in the azimuthal
component around the vortex point. A mathematical model for the interaction of
such vortex points in the plane was derived by Helmholtz (1858) and Kirchhoff
(1876). In bounded domains the equations for $n$ point vortices are derived
using the Kirchhoff-Routh functions; see \cite{Li41}. The first reference to
the model of vortices on a sphere can be found in \cite{Eb}. Recent
applications of the point vortices model include Bose-Einstein condensates and
semiconductors; see \cite{Ba17} and references therein. A relatively recent
exposition of topics related to the $n$-vortex problem can be found in
\cite{Ne01}, and the references therein.

The study of equal point vortices located at the vertices of a polygon, which
rotates around its center, originates in the work of Lord Kelvin (1878) and
Thomson (1883). They determined the linear stability properties of this
configuration in the plane. The linear stability of the polygonal relative
equilibrium for the case of a disk was studied much later \cite{Ha31}. The
nonlinear stability properties of the polygon in the plane were determined in
\cite{Ca99}, and for the case of the sphere in \cite{Bo00}. When the stability
of the $n$-polygon of vortices with a central vortex changes, then there is a
bifurcation of relative equilibria consisting of nested polygons in the plane
\cite{GaIz11}, and on the sphere \cite{Li01}.

A \textit{choreography} of the $n$-vortex problem is a solution where the $n$
vortices follow the same path. The principal aim of our article is to
determine such solutions in a systematic fashion, by making use of robust
boundary value continuation techniques in the presence of symmetries. The term
\textit{choreography} was adopted for the case of the $n$-body problem of
celestial mechanics after the work of Sim\'{o} \cite{Si00}. The first
non-circular choreography was discovered numerically for $3$ bodies in
\cite{Mo93}, and its existence was verified analytically in \cite{ChMo00},
using the direct method of calculus of variations. The results in
\cite{ChMo00} mark the beginning of the development of variational methods for
this purpose, where the existence of choreographies is associated with the
problem of finding critical points of the classical action of Newton's
equations of motion. However, in the case of the $n$-vortex problem it is
more difficult to establish the existence of choreographies using the direct
method of calculus of variations. The reason is that the action is strongly
indefinite and integrable at orbits with collisions, since $\ln r$ is
integrable at $r=0$.

Vortex choreographies can be constructed explicitly for $3$ vortices in the
plane, using the fact that the system is integrable in that case. Similarly,
choreographies for $4$ vortices in the plane have been found in \cite{Bo04}
and on the sphere in \cite{Bo05}, using the fact that the system can be
reduced to a system with few degrees of freedom. While the $3$-vortex problem
in the plane is integrable, the $n$-vortex problem for $n\geq4$ is only
integrable in special cases \cite{Ar82}. The absence of integrals of motion
makes it difficult to find explicit choreographies for $n\geq4$. For the case
of a disk, even the $3$-vortex problem is not integrable. Surprisingly,
choreographies have been constructed recently for $n$-vortices in general
bounded domains, using blow-up techniques. They are located close to a
stagnation point of a vortex in the domain \cite{Ba16} and close to the
boundary of the domain \cite{Ba18}.

Choreographies for the $n$-body problem exist along Lyapunov families that
arise from the $n$-polygon of bodies in a rotating frame \cite{CaDoGa,Ch09}.
Related results on choreographies in the discrete nonlinear Schr\"odinger
equations can be found in \cite{CaDoGaPa}. These results depend only on the
symmetries of the equations, and can be extended to the case of the $n$-vortex
problem in radially symmetric domains. The global existence of Lyapunov
families that arise from the $n$-polygon of vortices is established in the
plane in \cite{GaIz12}, and on the sphere in \cite{Ga18}. In our current paper
we determine choreographies along Lyapunov families of the $n$-vortex problem
in the plane, in a disk, and on a sphere. The families emanate from the
polygonal equilibrium of the vortices, with starting frequencies $\nu_{k}$
that are equal to the normal modes of oscillation of the equilibrium, and they
constitute continuous families in the space of renormalized periodic
functions. The \emph{global} property means that the Sobolev norm or the
period of the orbits along the family tends to infinity, or that the family
ends in a collision or at another equilibrium.

Specifically, let $\omega$ be the frequency of rotation of the polygonal
configuration in the plane, in a disk, or on a sphere, and let $\nu$ be the
frequency of the Lyapunov orbit in the rotating frame. A solution of the
$n$-vortex problem in a radially symmetric domain is given by $q_{j}%
(t)=e^{i\omega t}x_{j}(\nu t)$, where $x_{j}$ is a $2\pi$-periodic
renormalized Lyapunov orbit satisfying the symmetries
\begin{equation}
x_{j}(t)=e^{ij\zeta}x_{n}(t+jk\zeta),\qquad\zeta=2\pi/n~, \label{PS}%
\end{equation}
for some $k$, with $1\leq k\leq n$. A Lyapunov orbit is $\ell:m$ resonant when
$\ell$ and $m$ are relatively prime such that
\[
\frac{\omega}{\nu}=\frac{\ell}{m}~,\qquad k\ell-m\in n\mathbb{Z}~.
\]
These conditions on the frequency are satisfied for a dense set of parameter
values. An $\ell:m$ resonant Lyapunov orbit is a choreography in the inertial
frame and each of the integers $k$, $\ell$, and $m$ is of distinct
significance for the properties of the choreographies: the choreography has
winding number $\ell$ around a center, is symmetric with respect to rotations
by $2\pi/m$, and the $n$ vortices form groups of $d$-polygons, where $d$ is
the greatest common divisor of $k$ and $n$. For the numerical determination of
the choreographies we use boundary value continuation methods, as implemented
in the most recent version \cite{EJD07} of the AUTO software \cite{EJD81}.

The $n$-vortex problem is not homogeneous in the case of a disk and in the
case of a sphere. In these cases we can determine continuous families of
choreographies that correspond to $\ell:m$-resonant Lyapunov orbits, and that
are not related by a homogeneous scaling. This statement does not hold for the
case of vortices in the plane. We also consider partial choreographies, in
which an additional vortex of variable circulation is located at the center of
the polygon. The results of our paper can be extended to the case of $n$
vortices rotating in any radially symmetric manifold. We have also observed
choreographies where the vortices are arranged along a line. This phenomenon
may be related to the fact that a row of evenly arranged vortices of equal
vorticity is stable, according to \cite{Ar95} and \cite{Ar02}.

In Section~\ref{sec:plane} we recall how a dense set of orbits along the
Lyapunov families corresponds to choreographies. In Section~\ref{sec:numerics}
we describe in some detail the numerical continuation procedure used to
determine the periodic solution families that arise from the polygonal
relative equilibrium in the plane. Later sections include complementary
details on the specific numerical procedure used there. In
Section~\ref{sec:disk} we consider the $n$-vortex problem in a disk.
Section~\ref{sec:center} concerns the $\left( n+1\right) $-vortex problem,
with $n$ vortices of unit circulation, and with an additional vortex of
circulation $\mu$ at the center of the polygon. Section~\ref{sec:sphere}
concerns choreographies along Lyapunov families for the case of the $n$-vortex
problem on a sphere, for which we present somewhat more extensive numerical
results.

\section{The $n$-polygon of vortices in the plane}

\label{sec:plane}

Let $q_{j}(t)\in\mathbb{C}$ be the position of the $j$th vortex in the plane.
Assume that the $n$ vortices have equal circulation for $j=1,\ldots,n$. The
Hamiltonian of $n$ vortices in the plane is
\[
H(q)=-\frac{1}{2}\sum_{1\leq k<j\leq n}\ln\left\vert q_{j}-q_{k}\right\vert
^{2}~,
\]
and the symplectic form is $\omega=\sum_{j=1}^{n}dx_{j}\wedge dy_{j}$. Since
the Hamiltonian is invariant under the group of transformations $SO(2)$, which
acts as $e^{i\theta}u$, the equations also have the conserved quantity
\begin{equation}
G(q)=\frac{1}{2}\sum_{j=1}^{n}\left\vert q_{j}\right\vert ^{2}.
\end{equation}
The equations for the vortices in rotating coordinates, $q_{j}(t)=e^{i\omega
t}u_{j}(t)$, are then given by
\begin{equation}
i\dot{u}_{j}=\omega u_{j}-\sum_{k=1~(k\neq j)}^{n}\frac{u_{j}-u_{k}%
}{\left\vert u_{j}-u_{k}\right\vert ^{2}}~,\qquad u_{j}(t)\in\mathbb{C}%
~\text{.} \label{NE}%
\end{equation}
The $n$-polygon of vortices $u_{j}=e^{ij\zeta}$, $\zeta=2\pi/n$ , for
$j=1,\ldots,n$, is an equilibrium \cite{GaIz12} when
\begin{equation}
\omega=\frac{n-1}{2}~.
\end{equation}
In \cite{GaIz12}, equivariant degree theory is used to prove that the
polygonal relative equilibrium has a global family of periodic solutions of
the form $u_{j}(t)=x_{j}(\nu t)$, with symmetries (\ref{PS}), for each
positive frequency
\[
\nu_{k}=\left(  s_{k}\left(  2\omega-s_{k}\right)  \right)  ^{1/2}~,
\]
with $k=1,\ldots,n-1$, where $s_{k}=k(n-k)/2$.

Due to rotational invariance there is always a zero frequency, $\nu_{0}=0$.
Since the frequencies have $1:1$ resonances ($\nu_{k}=\nu_{n-k}$),
straightforward application of the Lyapunov center theorem is not possible,
except for some cases with  $k=n/2$ for $n=2,4,6$, see \cite{Ca14}. For
example, for $n=6$ we have
\[
\nu_{1}=\nu_{5}=\frac{5}{2}~,\qquad\nu_{2}=\nu_{4}=2~,\qquad\nu_{3}=\frac
{3}{2}~.
\]
According to \cite{CaDoGa}, we define a Lyapunov orbit as being $\ell:m$
\textit{resonant}~ if its frequency satisfies the relation%
\[
\frac{\omega}{\nu}=\frac{\ell}{m}~\text{,}%
\]
where $\ell$ and $m$ are relatively prime such that $k\ell-m\in n\mathbb{Z}$~.
The explicit period of an $\ell:m$ resonant Lyapunov orbit is
\[
T_{\ell:m}~=~\frac{2\pi}{\nu}~=~\frac{4\pi}{n-1}~\frac{\ell}{m}~.
\]

\begin{theorem}
\label{proposition} In the inertial frame an $\ell:m$ resonant Lyapunov orbit
is a choreography, satisfying
\[
q_{j}(t)=q_{n}(t+j\tilde{k}\zeta)~\text{,}%
\]
where $\tilde{k}=k-(k\ell-m)\ell^{\ast}$, with $\ell^{\ast}$ the $m$-modular
inverse of $\ell$. The period of the choreography is $m~T_{\ell:m}$. The
choreography is symmetric with respect to rotations by the angle $2\pi/m$, and
it winds around a center $\ell$ times.
\end{theorem}

\noindent For a proof of this result see \cite{CaDoGa} and \cite{Ga18}.

\section{Numerical continuation of Lyapunov families}

\label{sec:numerics}

For the numerical continuation of the Lyapunov families it is necessary to
take the rotational symmetries into account, \textit{i.e.}, to continue the
Lyapunov orbits numerically, we use the augmented equations
\begin{equation}
\dot{u}=(-i\omega+\lambda_{1})\nabla G-(i-\lambda_{2})\nabla H~. \label{EqA}
\end{equation}
The solutions of Equation~(\ref{EqA}) are solutions of the original
equations of motion when the {\it unfolding parameters} $\lambda_{1}$ and 
$\lambda_{2}$ are zero. 
The converse of this statement is also true, as long as the fields
$\nabla G$ and $\nabla H$ are linearly independent. 
The unfolding parameters and the corresponding unfolding terms in
Equation~\ref{EqA} are needed to regularize the continuation of periodic
solutions in this conservative system with rotational invariance; 
see \cite{DoVa03} for a general treatment of continuation of periodic 
orbits in conservative systems.

Below we provide key details on our computational approach, so that it can
also be applied to related problems, or be implemented in other computational
environments. Moreover, a directory with python scripts that make the AUTO
software carry out all computations reported in this paper, as well as much
more extensive computations, will be made freely available.

Due to the symmetries, most eigenvalues corresponding to the polygonal
equilibrium have multiplicity $2$. In order to compute emanating families of
periodic solutions, while avoiding lengthy mathematical derivations, we use a
simple perturbation approach. For the stationary solutions we use the
equations
\begin{equation}
\left(  \lambda_{1}-i\right)  \omega u_{j}+i\sum_{k\neq j}\kappa_{k}%
\frac{u_{j}-u_{k}}{\left\vert u_{j}-u_{k}\right\vert ^{2}}=0,\quad j=1,\ldots n~,
\label{PL1}%
\end{equation}
where $\omega=(n-1)/2$, and where $\lambda_{1}$ is an unfolding parameter 
that is used to regularize the initial continuation of stationary solutions.
Note that for stationary solutions only one such unfolding parameter is needed,
while two unfolding parameters are required for the continuation of periodic
solutions.
Also note that we have introduced parameters $\kappa_{j}$, $j=1,...,n$,
where $\kappa_{j}=1$ for $j=1,...,n-1$, but where we treat $\kappa_{n}$ as 
a perturbation parameter that is allowed to take on values different from $1$.
More specifically, the parameter $\kappa_{n}$ is used to perturb the 
circulation of the $n$-th vortex temporarily in the first few steps of the
algorithm that leads to choreographies.
We also add a simple constraint on the location of the $n$th vortex, which
removes the rotational invariance, namely,
\begin{equation}
\operatorname{Im}u_{n}=0~.\label{PL2}%
\end{equation}
Considered as real equations with real variables, Equations~(\ref{PL1}) and
(\ref{PL2}) together define a system of $2n+1$ equations. Conceptionally we
may consider the $u_{j}$ and $\lambda_{1}$ as $2n+1$ real variables, and
$\kappa_{n}$ as the continuation parameter, even though the continuation
algorithm does not make this distinction. The perturbation procedure carries
out this continuation, namely until $\kappa_{n}$ reaches a target value
different from $1$, for which we have used, for example, $\kappa_{n}=1.2$. For
the case $n=5$ that we consider in this section, the eigenvalues of the
linearized equations for the unperturbed equations are $\pm i\nu_{1}=\pm2i$
and $\pm i\nu_{2}=\pm i\sqrt{3}$, where each of these conjugate, purely
imaginary pairs has algebraic multiplicity $2$. In addition there is a zero
eigenvalue of multiplicity $2$. By contrast, the eigenvalues of the perturbed
equations, with $\kappa_{n}=1.2$, are found to be $\pm2i$, $\pm1.99717i$,
$\pm1.74814i$, and $\pm1.72256i$, as well as two real eigenvalues of opposite
sign that are very close to zero. These four purely imaginary eigenvalues are
simple, and each gives rise to a family of periodic solutions.

Each bifurcating family of periodic orbits is computed in three stages. 
The first stage is to follow the Lyapunov family starting from the 
perturbed equilibrium (with $\kappa_{n}$ close to, but different from $1$),
and until the amplitude of the periodic orbit reaches a small target value.
Here the ``amplitude'' is defined as 
\begin{equation}
A=\sum_{j=1}^{n}\int_{0}^{1}|u_{j}(t)-u_{j}^{0}|^{2}~dt~,\label{AMP}%
\end{equation}
where $u_{j}^{0}$ denotes the equilibrium position of the $j$th component.
Thus the amplitude is zero at the perturbed equilibrium.
In the second stage this
small amplitude orbit is followed keeping the amplitude fixed, while allowing
$\kappa_{n}$ to vary until it returns to the value $1$. This nontrivial
periodic orbit is then a solution of the unperturbed equations. In the third
stage this orbit is followed, keeping $\kappa_{n}$ fixed at $1$ and allowing
the amplitude to vary again, thereby generating the desired solution family of
the unperturbed problem.

Each of the three stages referred to above utilizes a boundary value
formulation for the continuation of periodic solutions. The differential
equation is here written as
\begin{equation}
\dot{u}_{j}=T(\lambda_{1}-i)\omega u_{j}-T(\lambda_{2}-i)\sum_{k\neq j}%
\kappa_{k}\frac{u_{j}-u_{k}}{\left\vert u_{j}-u_{k}\right\vert ^{2}},\quad
j=1,\ldots n~,\label{PL3}%
\end{equation}
where time has been rescaled to the unit interval $[0,1]$, so that the actual
period $T$ appears explicitly in the differential equations. The periodicity
boundary constraints are then given by
\begin{equation}
u_{j}(1)-u_{j}(0)=0~,\qquad j=1,\ldots,n~.\label{PBC}%
\end{equation}
The boundary value formulation also contains integral constraints. One of
these is the usual integral phase condition, here applied to the $n$th vortex
only, namely,
\begin{equation}
\operatorname{Re}\int_{0}^{1}u_{n}(t)~\dot{\tilde{u}}_{n}(t)~dt=0~,\label{PHS}%
\end{equation}
where $\dot{\tilde{u}}_{n}(t)$ represents the time derivative of a reference
solution, namely the preceding solution in the numerical continuation process.
Another integral constraint sets the average of the imaginary part of the
$n$th vortex to zero, namely,
\begin{equation}
\operatorname{Im}\int_{0}^{1}u_{n}(t)~dt=0~.\label{ROT}%
\end{equation}
This second integral constraint removes the rotational invariance of periodic
orbits. There are more general integral constraints for fixing the phase and
for removing invariances. However the ones listed above are simple, and
appropriate in the current context. Furthermore, one can keep track of the
\textquotedblleft amplitude\textquotedblright\ (\ref{AMP}) of the orbits.
More importantly, one can also choose to keep this amplitude fixed, 
by adding it as an integral constraint.

\begin{figure}[h]
\begin{center}
\resizebox{13.5cm}{!}{
\includegraphics{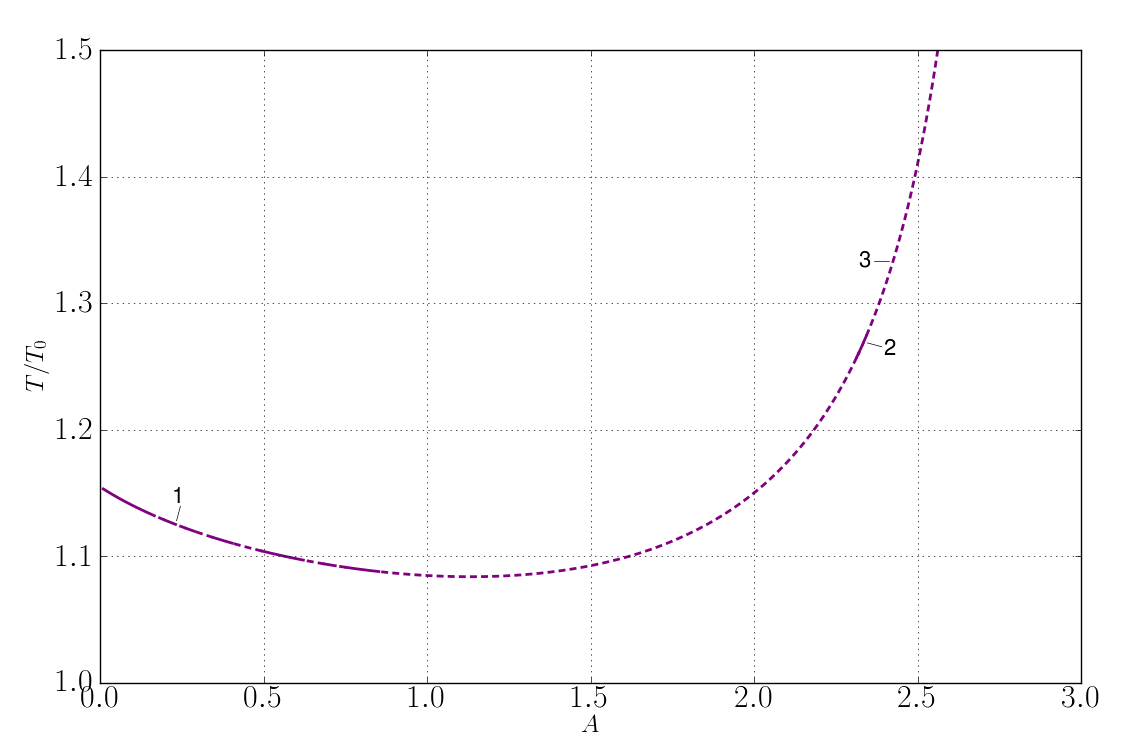}~}
\end{center}
\par
\vskip-.5cm\noindent\caption{ A representation of one of the four families of
periodic solutions for the case of 5 vortices in the plane, showing the
resonance ratio versus the amplitude parameter $A$. Solid curve sections
correspond to stable or almost stable solutions, and dashed sections
correspond to unstable solutions. The solutions labeled $1$, $2$, and $3$
correspond to the choreographies shown in Figure~\ref{fig02}. Solutions $1$
and $2$ are stable, with resonance ratio $9:8$, and $33:26$, repectively.
Solution $3 $ is unstable, with resonance ratio $4:3$. Further data on these
orbits can be found in Table~1. }%
\label{fig01}%
\end{figure}

\begin{figure}[ptb]
\begin{center}
\resizebox{18.0cm}{!}{
\includegraphics{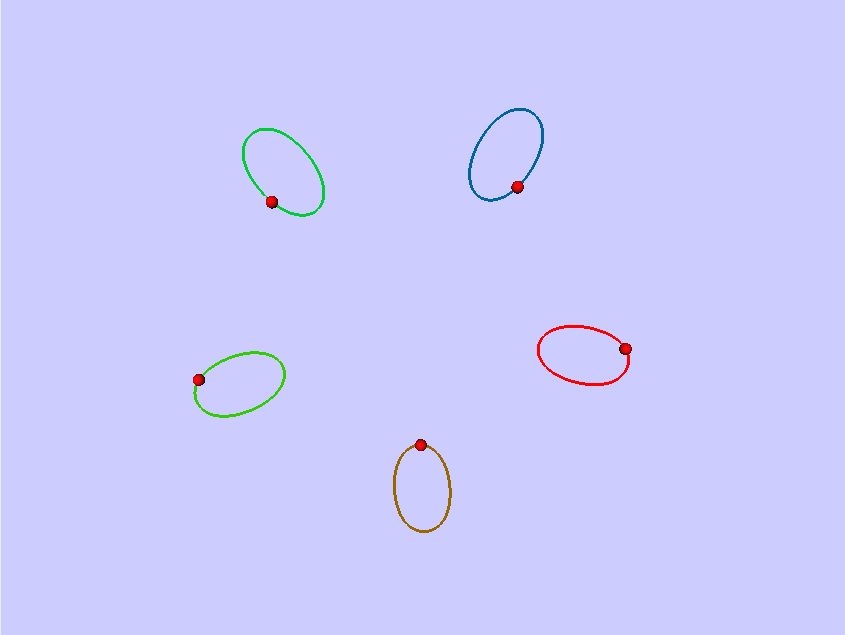}~
\includegraphics{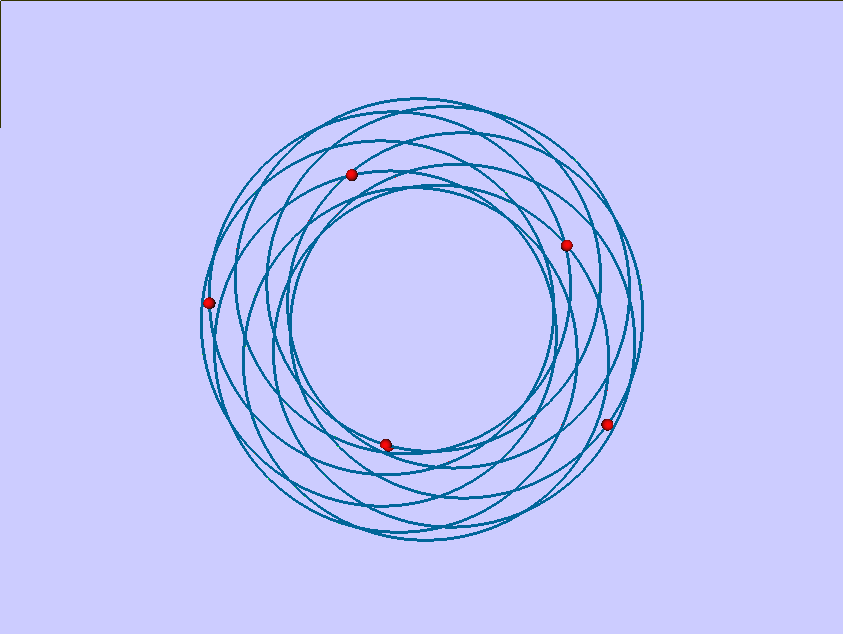} }
\end{center}
\par
\vskip-1.10cm\noindent
\par
\begin{center}
\resizebox{18.0cm}{!}{
\includegraphics{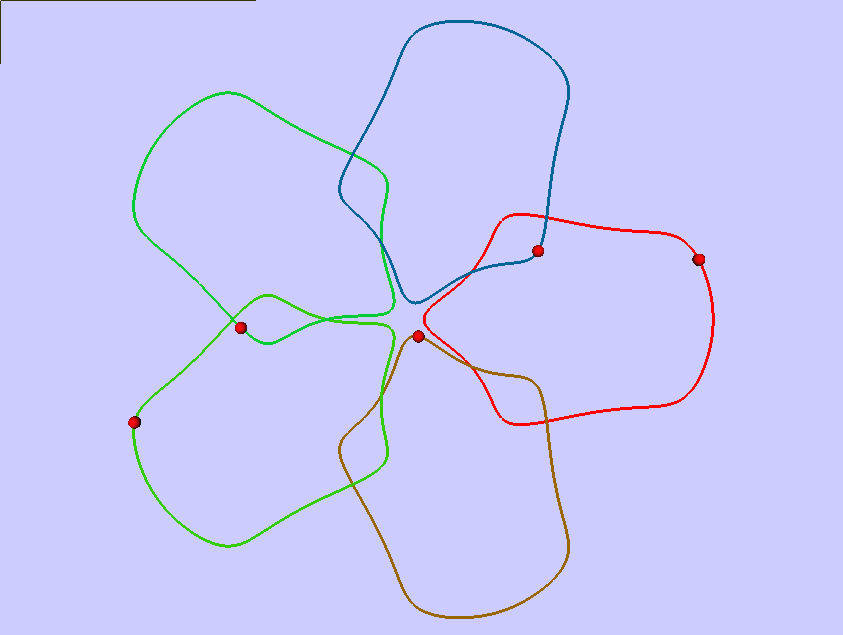}~
\includegraphics{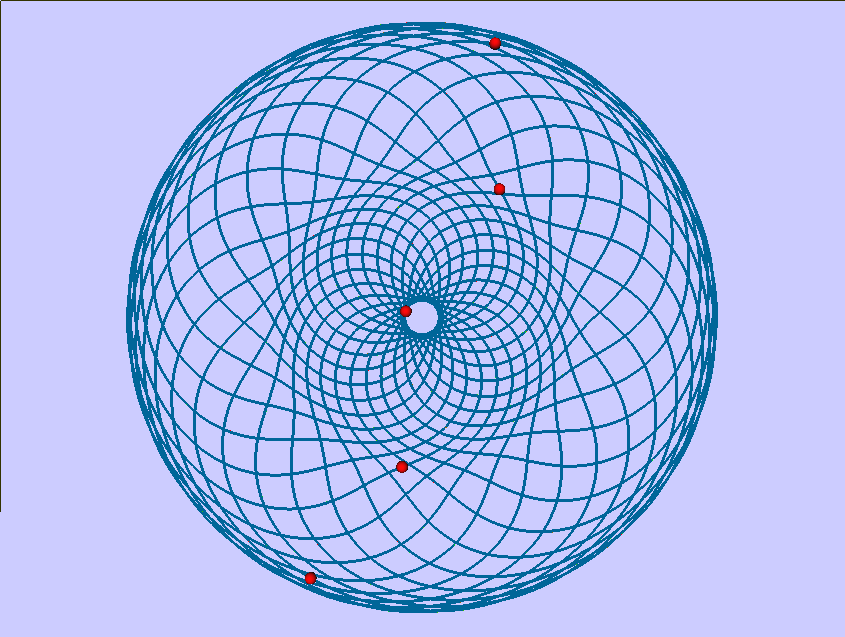} }
\end{center}
\par
\vskip-1.10cm\noindent
\par
\begin{center}
\resizebox{18.0cm}{!}{
\includegraphics{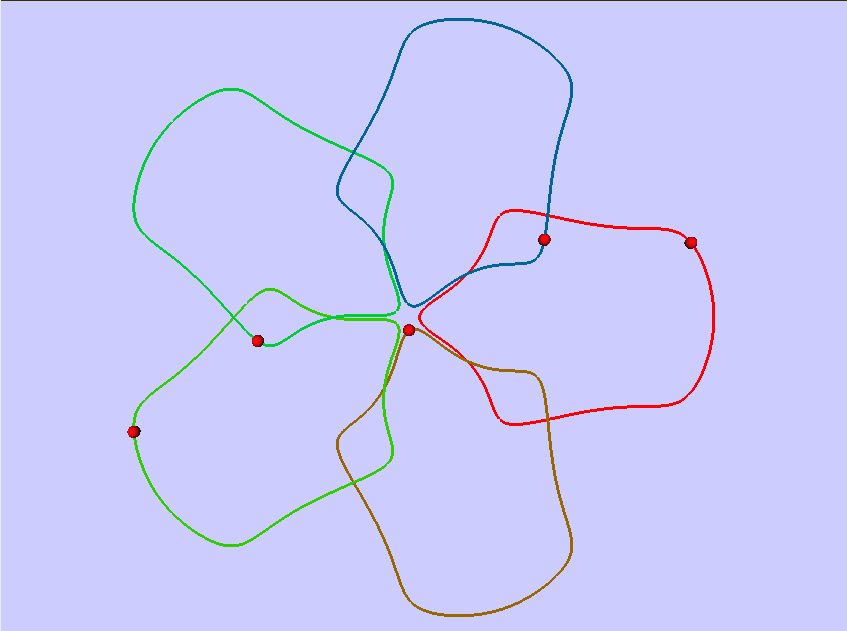}~
\includegraphics{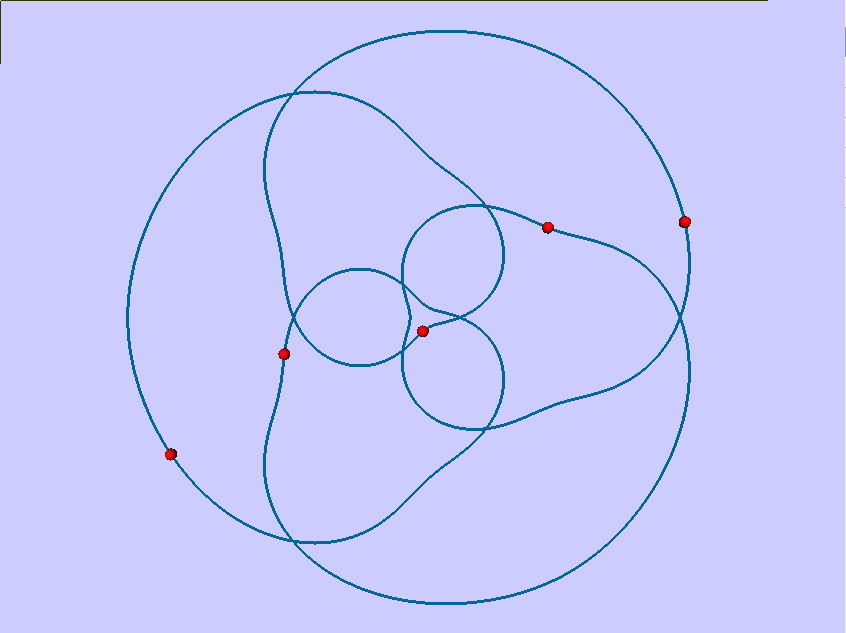} }
\end{center}
\par
\vskip-.5cm\noindent\caption{ The case of 5 vortices in the plane, for which a
bifurcation diagram is shown in Figure~\ref{fig01}. Numerical data can be
found in Table~1. The panels on the left show periodic solutions in the
rotating frame. The panels on the right show the corresponding periodic
solutions in the inertial frame, where they correspond to choreographies. }%
\label{fig02}%
\end{figure}

As outlined above, the determination of each bifurcating family of periodic
solutions is done in three stages. In the first stage, with perturbed
circulation $\kappa_{n}$, the free scalar parameters are the amplitude $A$,
the period $T$, and the unfolding parameters $\lambda_{1}$ and $\lambda_{2}$.
In the second stage, where the perturbation is being undone, the free
parameters are the circulation of the $n$th vortex, namely, $\kappa_{n}$, the
period $T$, and $\lambda_{1}$ and $\lambda_{2}$, while the amplitude $A$
remains fixed at a small, nonzero value. In the third stage, where the actual
family of interest is computed, the free parameters are again those used in
stage~1, namely, $A$, $T$, $\lambda_{1}$, and $\lambda_{2}$.
Figure~\ref{fig01} shows a representation of one of the four periodic solution
families that arise from the polygonal equilibrium for the case of $5$
vortices in the plane. More specifically, this family is one of the two
families that arise from the complex conjugate purely imaginary eigenvalue
$\pm i\sqrt{3}$, which has algebraic multiplicity $2$. Three representative
periodic solutions along this family can be seen in Figure~\ref{fig02}; each
shown in the rotating frame, as well as in the inertial frame, where they
correspond to choreographies.


\section{The vortex polygon in a disk}

\label{sec:disk}

Let $B_{R}\subset\mathbb{C}$ be the disk of radius $R$, and let $q_{j}(t)\in
B_{R}$ be the position of the $j$th vortex. Assume that the $n$ vortices have
equal circulation for $j\in\{1,\ldots,n\}$. According to \cite{Ku07}, the
Hamiltonian of $n$ vortices in the disk $B_{R}$ is given by
\[
H(q)=-\frac{1}{2}\sum_{k<j}\ln\left\vert q_{j}-q_{k}\right\vert ^{2}+\frac
{1}{4}\sum_{j,k=1}^{n}\ln\left\vert R^{2}-q_{j}\bar{q}_{k}\right\vert ^{2}~.
\]
The symplectic form and the conserved quantity $G(q)$ are the same as in the
case of the plane. The equations of motion are given by
\[
i\dot{q}_{j}=-\sum_{k\neq j}\frac{q_{j}-q_{k}}{\left\vert q_{j}-q_{k}%
\right\vert ^{2}}+\sum_{k=1}^{n}\frac{q_{j}-\hat{q}_{k}}{\left\vert q_{j}%
-\hat{q}_{k}\right\vert ^{2}}~,\qquad q_{j}(t)\in B_{R}~,
\]
where the second sum represents the interaction with the boundary, with
\[
\hat{q}_{k}=R^{2}/\bar{q}_{k}~.
\]
Since $\hat{q}_{k}=e^{i\omega t}\hat{u}_{k}$, with $\hat{u}_{k}=R^{2}/\bar
{u}_{k}$, the equations of the $n$ vortices in rotating coordinates,
$q_{j}(t)=e^{i\omega t}u_{j}(t)$, are given by
\begin{equation}
i\dot{u}_{j}=\omega u_{j}-\sum_{k\neq j}\frac{u_{j}-u_{k}}{\left\vert
u_{j}-u_{k}\right\vert ^{2}}+\sum_{k=1}^{n}\frac{u_{j}-\hat{u}_{k}}{\left\vert
u_{j}-\hat{u}_{k}\right\vert ^{2}}~,\qquad u_{j}(t)\in B_{R}~.
\end{equation}
For a proper choice of $\omega$, there is always an equilibrium of the form
$u_{j}=R_{0}e^{ij\zeta}$ with $\zeta=2\pi/n$. Without loss of generality we
can assume, after a rescaling if necessary, that $R_{0}=1$. We use the radius
$R$ of the disk as the problem parameter. The frequency of rotation of the
polygonal equilibrium $u_{j}=e^{ij\zeta}$, with $R>1$ \cite{Ha31}, is given
by
\begin{equation}
\omega=\frac{n}{1-R^{-2n}}-\frac{n+1}{2}~. \label{omega-disk}%
\end{equation}

The theorem in \cite{GaIz12} can also be applied to the polygonal relative
equilibrium in the disk. It implies that the polygon has a global family of
periodic solutions of the form $u_{j}(t)=x_{j}(\nu t)$ with symmetries
(\ref{PS}) for each normal mode of oscillation. These normal modes are the
purely imaginary eigenvalues of the linearization, which correspond to the
linearly stable eigenvalues analyzed in \cite{Ha31}. Theorem~\ref{proposition}
and its numerical implementation in the plane are applicable in the disk as
well, \textit{i.e.}, the $\ell:m$ resonant Lyapunov orbits in the disk, with
period $T_{\ell:m}=2\pi\frac{\ell}{\omega m}$, correspond to choreographies in
the inertial frame.

The augmented equations for the numerical determination of the families of
periodic solutions for the case of a disk are given by
\[
\dot{u}_{j}=T(\lambda_{1}-i)~\omega u_{j}+T(\lambda_{2}-i)\left(  \sum_{k\neq
j}\frac{u_{j}-u_{k}}{\left\vert u_{j}-u_{k}\right\vert ^{2}}-\sum_{k=1}%
^{n}\frac{u_{j}-\hat{u}_{k}}{\left\vert u_{j}-\hat{u}_{k}\right\vert ^{2}%
}\right)  ~,
\]
which include two unfolding terms, multiplied by their respective unfolding
parameters $\lambda_{1}$ and $\lambda_{2}$.

For the case of vortices in a disk, the numerical determination of the
families of periodic orbits that arise from the polygonal equilibrium is
easier than for the case of the plane. The reason is that the purely imaginary
eigenvalues are now simple, although they approach higher multiplicity when
the radius parameter $R$ gets large. We also observe that for small values of
$R$ there can be unstable eigenvalues. One of the several cases that we have
considered in detail is that of $5$ vortices in a disk. When $R=6$ then the
purely imaginary eigenvalues for the case $n=5$ are given by
$\{1.99999i,~1.86111i,~1.73111i,~1.72371i\}$. Thus the perturbation approach
used for the case of vortices in the plane is not needed now. Apart from that,
the formulation of the boundary value problem for computing the families of
periodic solutions is similar to that for the case of the plane. In
particular, the same periodicity boundary conditions and integral constraints
can be used.

The presence of a problem parameter, namely $R$, allows for additional
continuation sche\-mes; for example, where a solution quantity such as the
period $T$ or the resonance ratio $T/T_{0}$ is fixed, while $R$ is free to
vary. As a variation on this, we have found it useful here to implement a
continuation scheme where the $n$th vortex at time zero, \textit{i.e.},
$u_{n}(0)$, is constrained to be located on the $x$-axis and near the origin;
see, for example, the panels on the left in Figure~\ref{fig04}. To this end we
add the following two boundary conditions to the periodicity
Equations~(\ref{PBC}):
\begin{equation}
\operatorname{Im}u_{n}(0)~=~0 \quad\mbox{and}\quad\operatorname{Re}u_{n}(0)
~-~ r_{n} ~=~0~,\label{BC}%
\end{equation}
which also introduces a parameter $r_{n}$, which can be chosen to remain fixed
or be free to vary. The condition $\operatorname{Im}u_{n}(0)=0$ replaces the
integral phase condition~(\ref{PHS}), while the integral constraint
(\ref{ROT}) that inhibits rotation remains included. As a result, the
parameter $r_{n}$ then corresponds to the closest approach of the $n$th vortex
to the origen. By symmetry, the value of $r_{n}$ also corresponds to the
closest approach to the origen of the other vortices. Furthermore, the
integral constraint (\ref{AMP}) that keeps track of the amplitude of the
orbits, also remains included. Effectively, the continuation system then
includes $2n+2$ real boundary constraints and $2$ real integral constraints.
Since the system of differential equations corresponds to $2n$ real equations,
there must be $5$ free scalar parameters to compute a 1-dimensional continuum
of periodic orbits. Periodic orbits are first continued with $r_{n}$, $A$,
$T$, $\lambda_{1}$, and $\lambda_{2}$ as the five free parameters, namely until
$r_{n}$ reaches a desired target value. Subsequently this target orbit is
continued keeping $r_{n}$ fixed, while allowing the problem parameter $R$ to
vary. The five free continuation parameters in this step are therefore $R$,
$A$, $T$, $\lambda_{1}$, and $\lambda_{2}$.

\begin{figure}[h]
\begin{center}
\resizebox{15.5cm}{!}{
\includegraphics{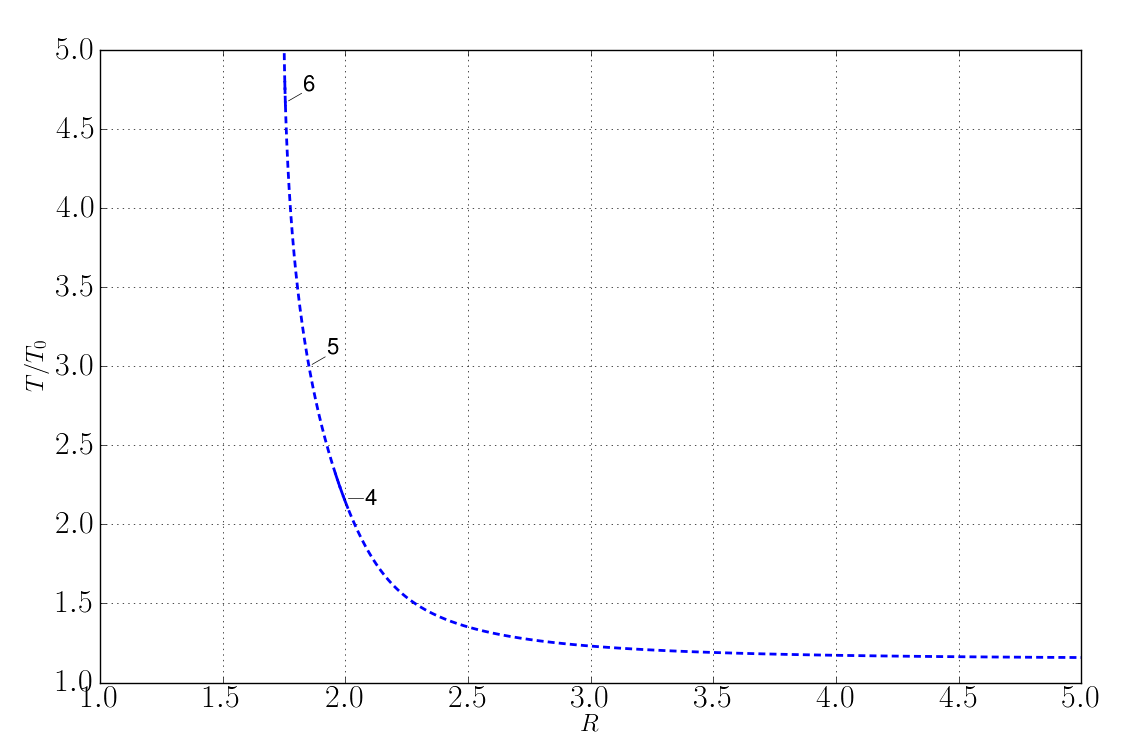}~}
\end{center}
\par
\vskip-.5cm\noindent\caption{ One of the four families of periodic solutions
that bifurcate from the polygonal equilibrium for the case of 5 vortices in a
disk. Solid curve sections represent stable or almost stable solutions. 
In particular, solutions~4 and 6 are seen to be stable, while solution~5 is 
unstable. The solutions labeled $4$, $5$, and $6$ correspond to choreographies, 
with the actual orbits shown in Figure~\ref{fig04}. As explained in the text, 
all orbits along this family share the property that their closest approach 
to the origin is equal to $r_{n}=0.2$. Further data on the labeled solutions 
can be found in Table~1. }%
\label{fig03}%
\end{figure}

\begin{figure}[ptb]
\begin{center}
\resizebox{18.0cm}{!}{
\includegraphics{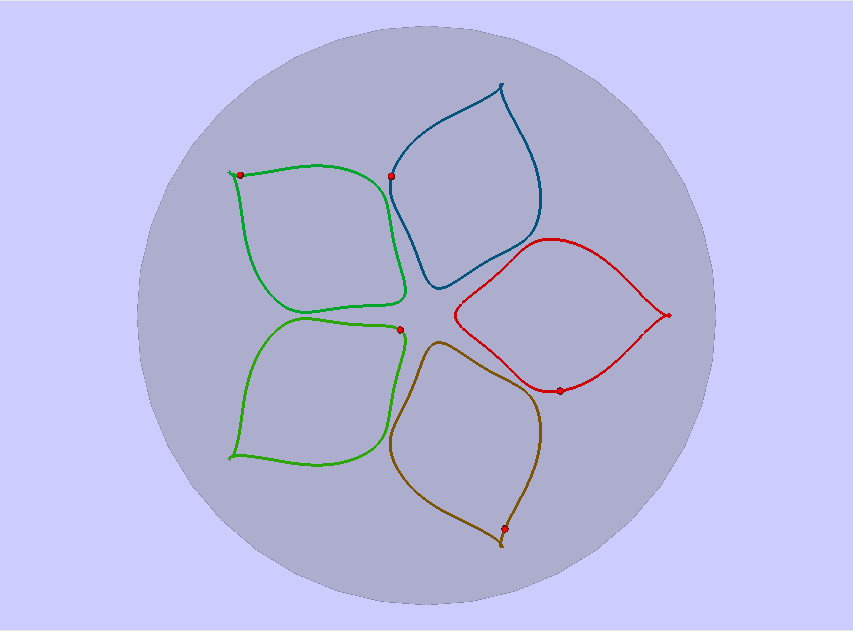}~
\includegraphics{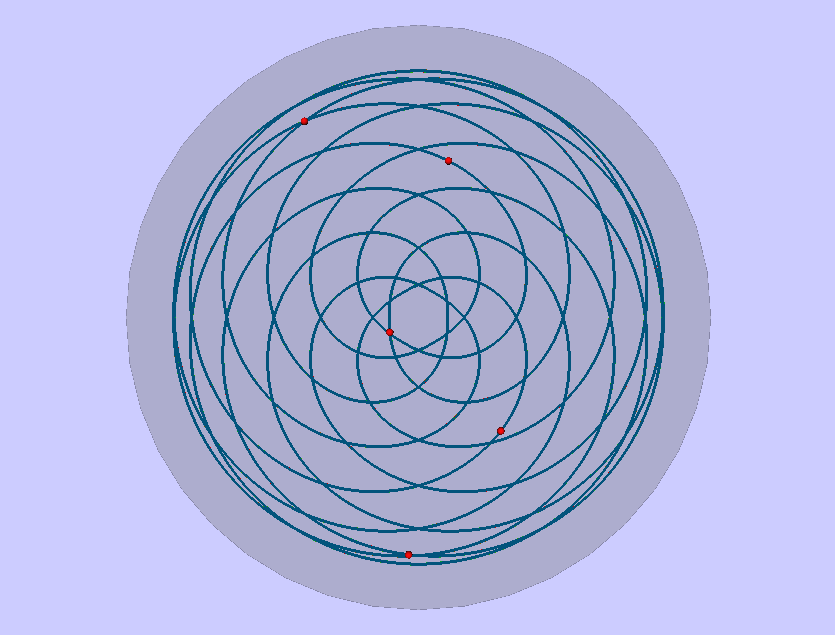} }
\end{center}
\par
\vskip-1.10cm\noindent
\par
\begin{center}
\resizebox{18.0cm}{!}{
\includegraphics{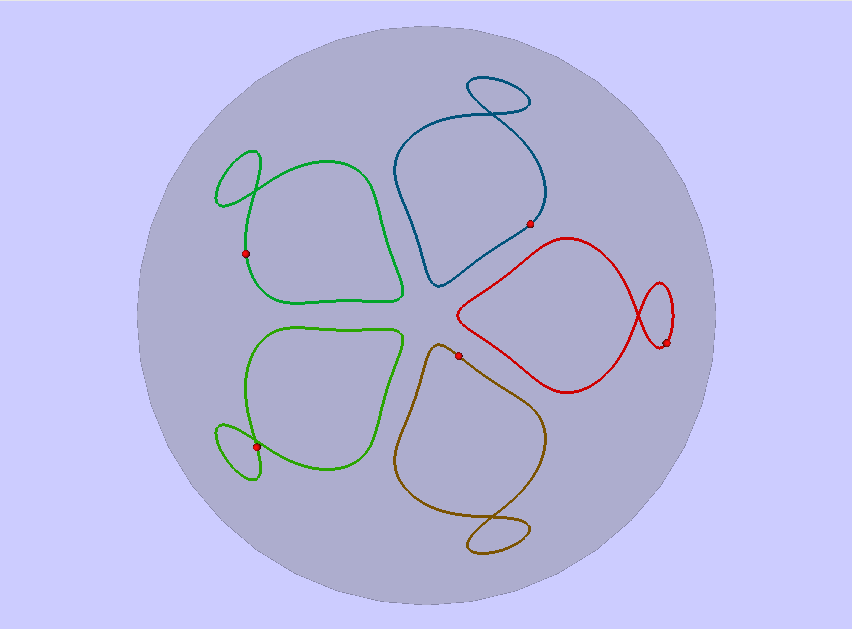}~
\includegraphics{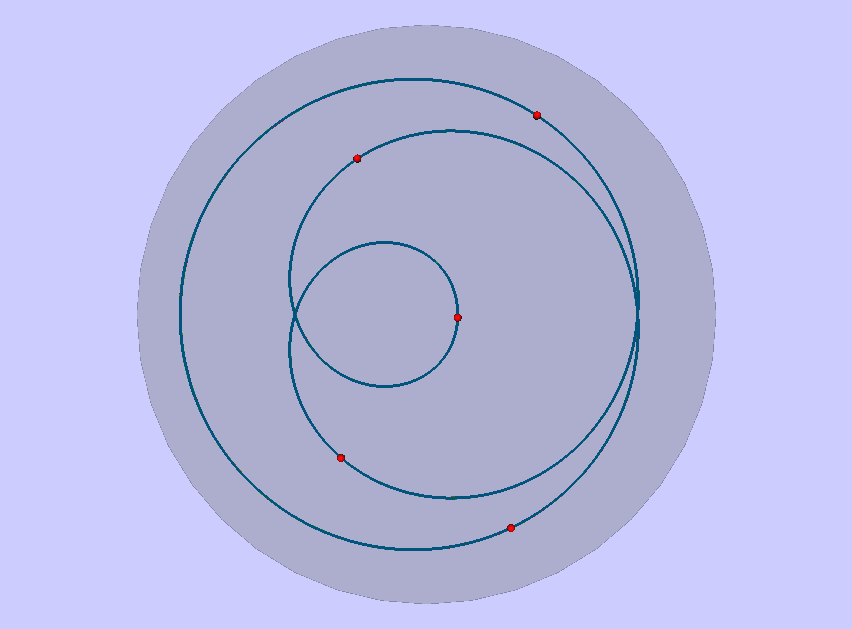} }
\end{center}
\par
\vskip-1.10cm\noindent
\par
\begin{center}
\resizebox{18.0cm}{!}{
\includegraphics{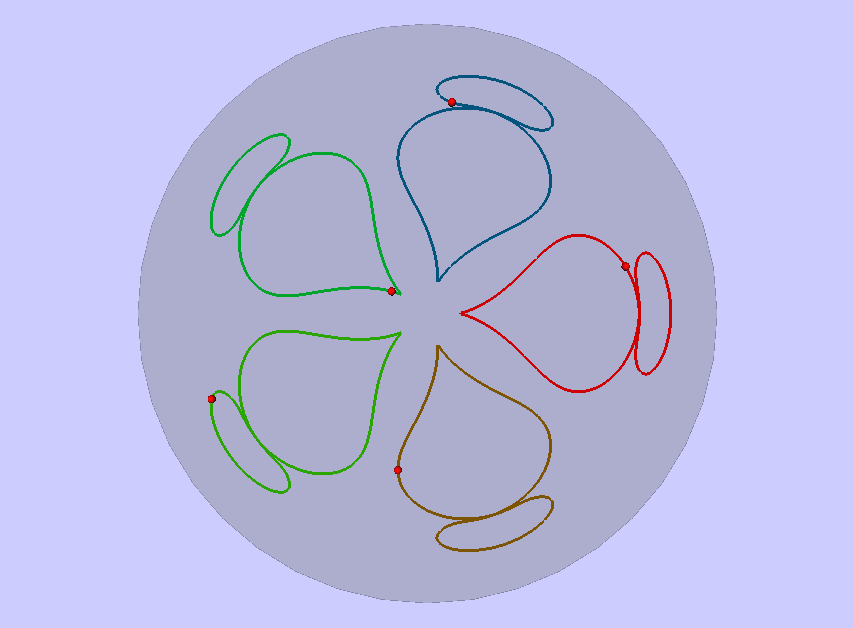}~
\includegraphics{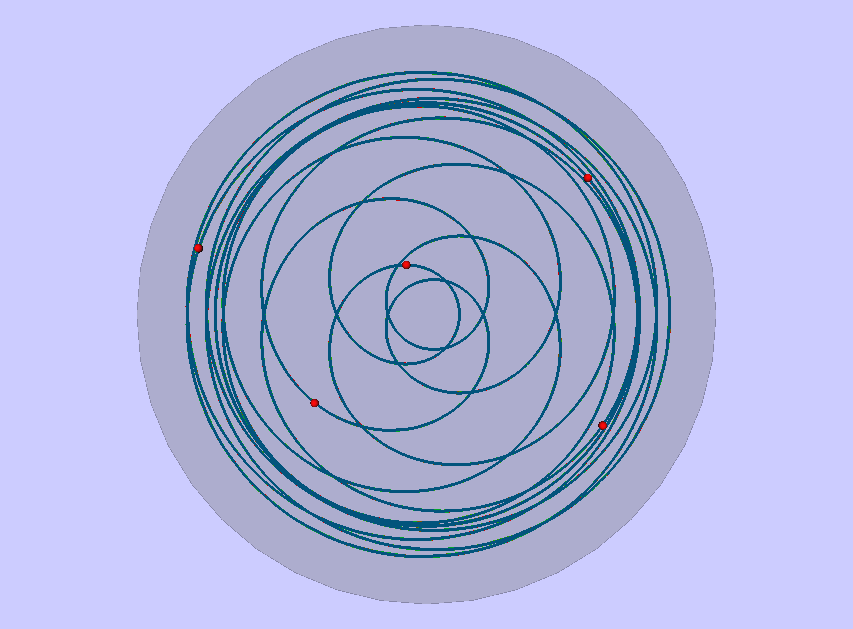} }
\end{center}
\par
\vskip-.5cm\noindent\caption{ The case of $n$ vortices in a disk, where $n=5$,
for which a bifurcation diagram is shown in Figure~\ref{fig03}. Numerical data
can be found in Table~1. The three rows of two panels correspond to the
solutions labeled $4$, $5$, and $6$, respectively, in Figure~\ref{fig03}. The
panels on the left show periodic solutions in the rotating frame. Here the
path of the $n$th vortex is symmetric with respect to the real axis, and its
closest approach to the origen occurs at time $t=0$, as a result of the
boundary and integral constraints described in the text. The panels on the
right show the corresponding periodic solutions in the inertial frame, where
they correspond to choreographies. }%
\label{fig04}%
\end{figure}

A representative bifurcation diagram for the case $n=5$, with $r_{n}=0.2$, is
shown in Figure~\ref{fig03}. The family that is shown there arises indirectly
from the purely imaginary eigenvalue $1.72371i$ of the polygonal equilibrium.
The family that emanates directly from this equilibrium is followed until it
reaches the target orbit having $r_{n}=0.2$. Subsequently the target orbit is
followed keeping $r_{n}$ fixed, while allowing the problem parameter $R$ to
vary. Also varying in this secondary continuation is the period $T$, while the
base period $T_{0}$ of the rotating frame is a direct function of $R$, namely
$T_{0}=2 \pi/\omega$, with $\omega=\omega(R)$ defined in
Equation~\ref{omega-disk}. Figure~\ref{fig03} provides a representation of the
solution family that is generated in the secondary continuation, by showing
the value of the resonance ratio $T/T_{0}$ versus the parameter $R$.

Along the solution path in Figure~\ref{fig03} there exists a countably
infinite number of resonant orbits that correspond to choreographies in the
inertial frame. Numerically we can in principle locate an unlimited number of
these. Three such choreographic solutions are indicated by the labeled points
in Figure~\ref{fig03}. The resonance ratio of these solutions is $13:6$ (label
4), $3:1$ (label 5), and $14:3$ (label 6), respectively. The actual orbits of
the three labeled solutions in Figure~\ref{fig03} are shown in
Figure~\ref{fig04}. The panels on the left in Figure~\ref{fig04} show these
resonant orbits in the rotating frame, and the panels on the right show the
corresponding orbits in the inertial frame, where they correspond to
choreographies. Observe that in Figure~\ref{fig04} the path of the $n$th
vortex, where $n=5$, is the one that is symmetric with respect to the real
axis. The closest approach to the origen occurs at time $t=0$, as a result of
the boundary and integral constraints described above.

Data on the solutions labeled 4, 5, and 6 in Figure~\ref{fig03} is given in
Table~1. In particular, solutions~4 and 6 are seen to be stable, while
solution~5 is unstable, as can also be inferred from Figure~\ref{fig03}. The
magnitude of the largest Floquet multiplier of the unstable solution~5 is
found to be $7.39340$. The particular manner by which these solutions were
obtained, as described above, represents one of various computational schemes
that we have used in the quest to locate stable solutions. As illustrated in
Figure~\ref{fig03}, regions of stability along solution families can be
relatively small, and therefore require a systematic approach for their
determination.

\section{The vortex polygon with a center vortex}

\label{sec:center}

Let $q_{j}(t)\in\mathbb{C}$ be the position of the $j$th vortex in the plane,
for $j=0,1,\ldots,n$. Assume the $n$ vortices with $j=1,\ldots,n$ have
circulation $1$, and the vortex with $j=0$ has circulation $\mu$. The
Hamiltonian is
\[
H(q)=-\frac{1}{2}\sum_{1\leq k<j\leq n}\ln\left\vert q_{j}-q_{k}\right\vert
^{2}-\frac{1}{2}\mu\sum_{k=1}^{n}\ln\left\vert q_{0}-q_{j}\right\vert ^{2}~,
\]
with symplectic form
\[
\omega=\mu\left( dx_{0}\wedge dy_{0}\right)  + \sum_{j=1}^{n}dx_{j}\wedge
dy_{j}~,
\]
and conserved quantity
\[
G(q)=\mu\frac{1}{2}\left\vert q_{0}\right\vert ^{2}+\frac{1}{2}\sum_{j=1}%
^{n}\left\vert q_{j}\right\vert ^{2}~.
\]
The equations of motion of the $n$ vortices in rotating coordinates,
$q_{j}(t)=e^{i\omega t}u_{j}(t)$, are given by
\begin{equation}
i\dot{u}_{j}=\omega u_{j}-\mu\frac{u_{j}-u_{0}}{\left\vert u_{j}%
-u_{0}\right\vert ^{2}}-\sum_{k\neq j}\frac{u_{j}-u_{k}}{\left\vert
u_{j}-u_{k}\right\vert ^{2}}~,\quad j=1,\ldots,n~,
\end{equation}%
\[
i\dot{u}_{0}=\omega u_{0}-\sum_{k=1}^{n}\frac{u_{0}-u_{k}}{\left\vert
u_{0}-u_{k}\right\vert ^{2}}~.
\]
The vortex ring with a center corresponds to the positions $u_{0}=0$, and
$u_{j}=e^{ij\zeta}$, $j=1,\ldots,n$. This configuration is an equilibrium
when
\begin{equation}
\omega=\mu+s_{1}=\mu+\frac{n-1}{2}~.\label{omega-center}%
\end{equation}
The configuration has a global family of periodic solutions of the form
$u_{j}(t)=x_{j}(\nu t)$, with symmetries (\ref{PS}) for each positive
frequency $\nu_{k}=\left(  s_{k}(2\omega-s_{k})\right)  ^{1/2}$, where
$s_{k}=k(n-k)/2$ \cite{GaIz12}. Theorem~\ref{proposition}, and the numerical
algorithms for $n$ equal vortices, are also applicable to the current case 
of $n+1$ vortices. Here an $\ell:m$ resonant Lyapunov orbit has period
\[
T_{\ell:m}=\frac{2\pi}{\omega}\frac{\ell}{m}~.
\]
For these resonant orbits, $n$ vortices move in a choreographic fashion in the
inertial frame, while the central vortex follows a separate closed curve.

The augmented system of equations used for numerical continuation of periodic
orbits is now
\begin{equation}
\dot{u}_{j}=T\left(  \lambda_{1}-i\right)  \omega u_{j}-T\left(  \lambda
_{2}-i\right)  \left(  \mu\frac{u_{j}-u_{0}}{\left\vert u_{j}-u_{0}\right\vert
^{2}}+\sum_{k\neq j}\frac{u_{j}-u_{k}}{\left\vert u_{j}-u_{k}\right\vert ^{2}%
}\right)  ~,\quad j=1,\ldots,n~,
\end{equation}%
\[
\dot{u}_{0}=T\left(  \lambda_{1}-i\right)  \omega u_{0}-T\left(  \lambda
_{2}-i\right)  \sum_{k=1}^{n}\frac{u_{0}-u_{k}}{\left\vert u_{0}%
-u_{k}\right\vert ^{2}}~.
\]
\begin{figure}[ptbh]
\begin{center}
\resizebox{13.5cm}{!}{
\includegraphics{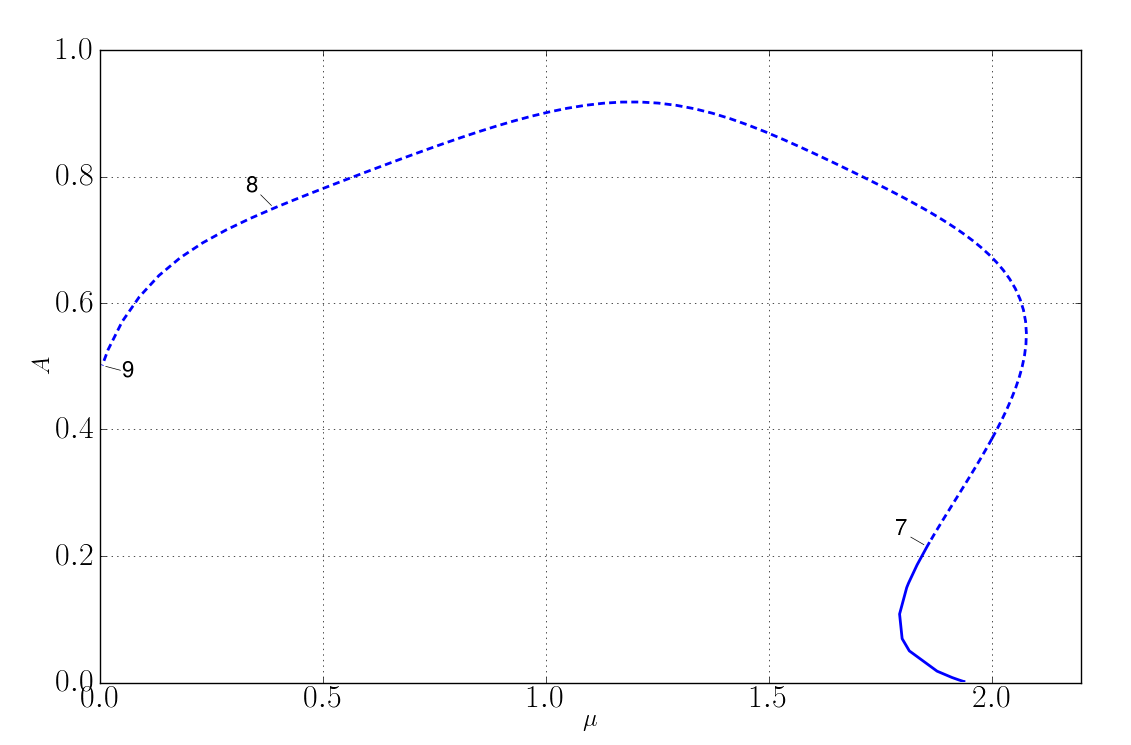}}
\end{center}
\par
\vskip-0.7cm\noindent\caption{ A continuum of choreographies of fixed
resonance $T/T_{0}=11:4$, with varying circulation $\mu$ of the central
vortex, for the case of 5 vortices with a 6th vortex at the center. The
solutions labeled $7$, $8$, and $9$ correspond to the partial choreographies
shown in Figure~\ref{fig06}. All have the same resonance ratio, namely $11:4$,
and only solution $7$ is stable. }%
\label{fig05}%
\end{figure}

As in Section~\ref{sec:plane} for vortices in the plane, we use a perturbation
approach to deal with the purely imaginary eigenvalues of multiplicity~2,
which also arise here. As in Section~\ref{sec:disk} for vortices in a disk,
the presence of a problem parameter allows for additional continuation
schemes. Here it was found useful to continue periodic orbits with fixed
resonance ratio $T/T_{0}$, allowing the circulation $\mu$ of the center vortex
to vary. To be more precise, via the perturbation approach we can compute all
families that arise from the purely imaginary eigenvalues. Along these
families we accurately detect selected resonances of interest; in particular
those that correspond to choreographies in the inertial frame. In follow-up
computations, such resonant orbits can be continued with fixed resonance ratio
$T/T_{0}$, allowing the circulation parameter $\mu$ to vary. With the
constraint that $T/T_{0}$ remain fixed, the standard periodicity boundary
conditions, and with three integral constraints, namely the integral phase
constraint applied to the $n$th vortex, the integral that sets the average
$y$-coordinate of the $n$th vortex to zero, and the integral that defines the
amplitude measure $A$, the full list of free continuation parameters is then
given by $\mu$, $A$, $T$, $\lambda_{1}$, and $\lambda_{2}$. Note that the
period $T_{0}$ of the rotating frame is taken to be a function of $\mu$,
namely $T_{0}=2 \pi/\omega$, with $\omega=\omega(\mu)$ defined in
Equation~\ref{omega-center}. In this way we can generate a \textit{continuum}
of choreographies, which is useful in the search for stable ones.

\begin{figure}[ptb]
\begin{center}
\resizebox{18.0cm}{!}{
\includegraphics{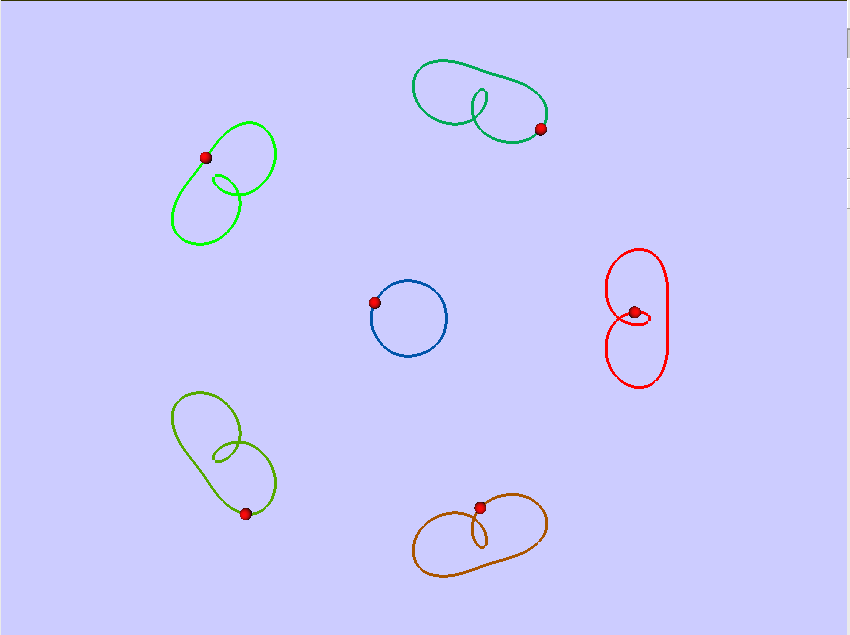}~
\includegraphics{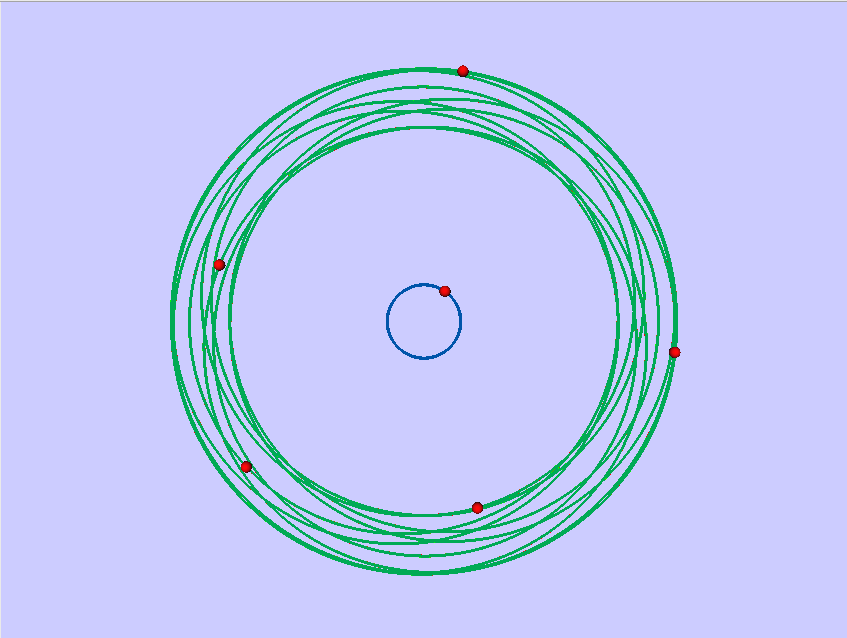} }
\end{center}
\par
\vskip-1.10cm\noindent
\par
\begin{center}
\resizebox{18.0cm}{!}{
\includegraphics{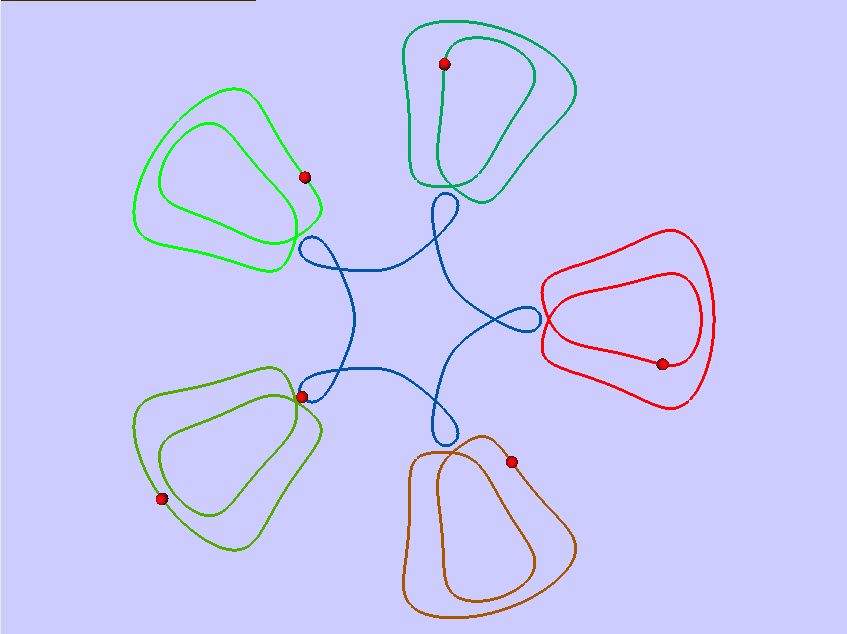}~
\includegraphics{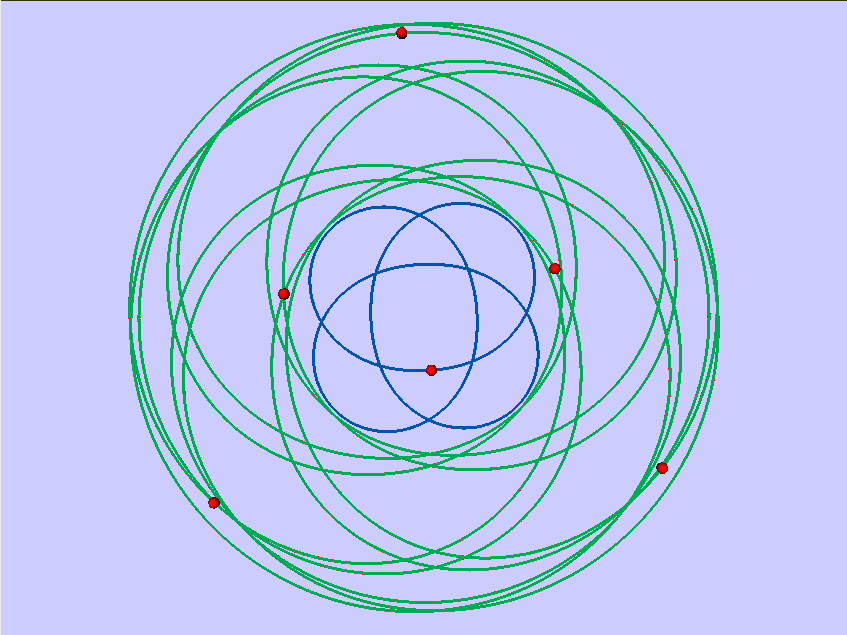} }
\end{center}
\par
\vskip-1.10cm\noindent
\par
\begin{center}
\resizebox{18.0cm}{!}{
\includegraphics{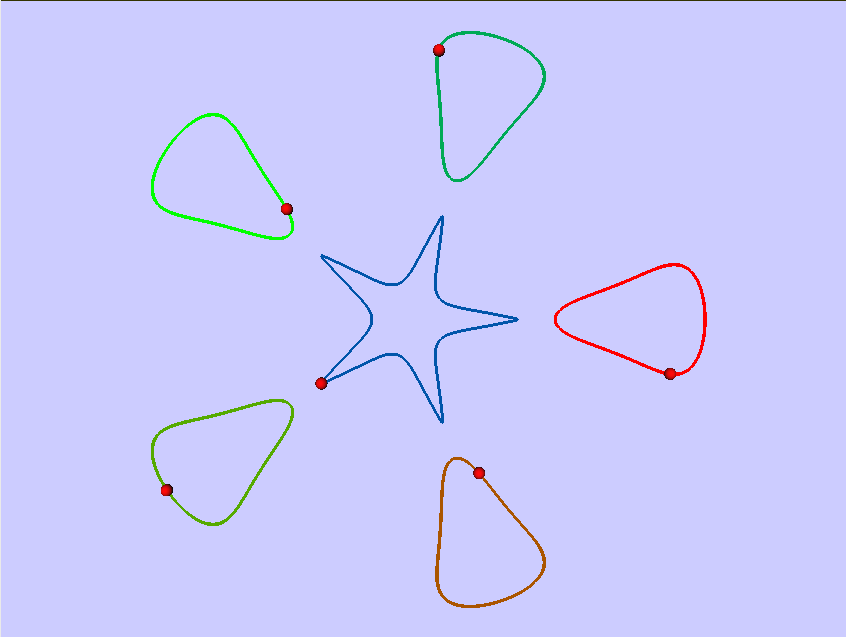}~
\includegraphics{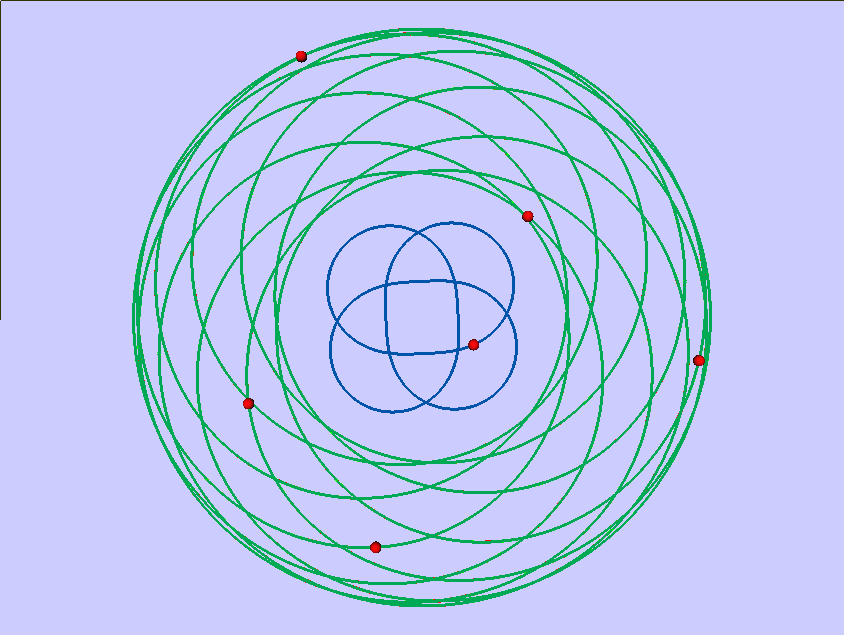} }
\end{center}
\par
\vskip-.5cm\noindent\caption{ The case of 5 vortices of unit circulation in
the plane, with a 6th vortex of variable circulation at the center. A
bifurcation diagram is shown in Figure~\ref{fig05}. Numerical data can be
found in Table~1. The panels on the left show periodic solutions in the
rotating frame. The panels on the right show the corresponding periodic
solutions in the inertial frame, where they are seen to correspond to partial
choreographies. }%
\label{fig06}%
\end{figure}

Figure~\ref{fig05} shows such a continuum, \textit{i.e.}, a family of
choreographies along which the resonance ratio remains constant, namely
$T/T_{0}=11:4$. A selection of three choreographies along this family is
indicated in Figure~\ref{fig05}, with the actual orbits shown in
Figure~\ref{fig06}, both in the rotating and in the inertial frame. The orbit
labeled 7 in Figure~\ref{fig05} is stable, or almost stable, with all
multipliers on the unit circle in the complex plane, to at least five decimal
digits accuracy. Orbit~8 is highly unstable, with two real multipliers outside
the unit circle, one having magnitude $1.3205~10^{4}$. Orbit~9 is for a very
small value of $\mu$, namely $\mu=0.005$, and is also unstable, with a complex
conjugate pair of multipliers of magnitude $38.1319$. \textit{A posteriori}
numerical integration shows that this periodic orbit is followed for a
relatively short time, after which the weak central vortex is ejected to the
outside.

\section{Vortices on a sphere}

\label{sec:sphere}

According to \cite{Ne01}, the motion of $n$ vortices on the sphere
$S^{2}=\{X\in\mathbb{R}^{3}:\left\Vert v\right\Vert =R\}$ is described by
vectors $v_{j}(t)\in S^{2}$ that satisfy%
\[
\dot{v}_{j}=\frac{1}{R}\sum_{k\neq j}\frac{v_{k}\times v_{j}}{\left\Vert
v_{j}-v_{k}\right\Vert ^{2}}~.
\]
Letting $v_{j}=(X_{j},Y_{j},Z_{j})$, the stereographic projection of the
sphere $S^{2}$ with $R=1$ onto the complex plane is%
\[
q_{j}=\frac{X_{j}+iY_{j}}{1-Z_{j}}\in\mathbb{C}~\text{.}%
\]
According to \cite{Mo13}, the Hamiltonian of the $n$ vortices on the sphere,
parameterized by the stereographic projection, is given by
\[
H(q) = -\frac{1}{2}\sum_{k\neq j}  \ln\frac{\left\vert q_{j}-q_{k}\right\vert
^{2}} {\left(  1+\left\vert q_{j}\right\vert ^{2}\right)  \left(  1+\left\vert
q_{k}\right\vert ^{2}\right)  }~,
\]
where the simplectic form is
\[
\omega= \sum_{j=1}^{n}  \frac{4}{\left( 1+\left\vert q_{j}\right\vert
^{2}\right) ^{2}}  ~dx_{j} \wedge dy_{j}~,  \qquad q_{j} = x_{j} + iy_{j}
\text{~.}%
\]
Here the conserved quantity is
\begin{equation}
G(q) = \sum_{j=1}^{n} \frac{\left\vert q_{j}\right\vert ^{2}} {1+\left\vert
q_{j} \right\vert ^{2}}~.
\end{equation}
The equations of the $n$ vortices on the sphere, parameterized by the
stereographic projection, are then given by
\begin{equation}
i\dot{q}_{j} = -\frac{1}{4} \left( 1+\left\vert q_{j}\right\vert ^{2}\right)
^{2} \sum_{k\neq j} \left( \frac{q_{j}-q_{k}}{\left\vert q_{j}-q_{k}%
\right\vert ^{2}} - \frac{q_{j}}{1+\left\vert q_{j}\right\vert ^{2}}\right)
\text{.}%
\end{equation}
In rotating coordinates, with $q_{j}(t)=e^{i\omega t}u_{j}(t)$, we have
\begin{equation}
i\dot{u}_{j}=\omega u_{j}-\frac{1}{4}\left(  1+\left\vert u_{j}\right\vert
^{2}\right)  ^{2}\sum_{k\neq j}\left(  \frac{u_{j}-u_{k}}{\left\vert
u_{j}-u_{k}\right\vert ^{2}}-\frac{u_{j}}{1+\left\vert u_{j}\right\vert ^{2}%
}\right)  \text{.} \label{Sp}%
\end{equation}
The $n$ vortices on the sphere have a polygonal equilibrium
\[
u_{j}=re^{ij\zeta}~,\qquad\zeta=2\pi/n~,\quad j=1,\ldots,n~,
\]
for each $r\in\lbrack0,1]$, with
\begin{equation}
\omega=\frac{n-1}{2}\left(  \frac{1-r^{4}}{4r^{2}}\right)
~.\label{omega-sphere}%
\end{equation}
A theorem in \cite{Ga18} establishes the existence of Lyapunov families of
periodic orbits that arise from polygonal relative equilibria on the sphere.
Specifically, each polygon has a global family of periodic solutions of the
form $u_{j}(t)=x_{j}(\nu t)$, with symmetries (\ref{PS}), for each normal mode
of oscillation. These normal modes of oscillation are the purely imaginary
complex eigenvalues of the linearization.

Theorem~\ref{proposition} and the numerical continuation algorithms for the
case of the plane are applicable to the sphere as well; \textit{i.e.}, 
there are $\ell:m$ resonant Lyapunov orbits with period 
$T_{\ell:m}=\frac{2\pi}{\omega}\frac{\ell}{m}$ 
that correspond to choreographies in the inertial frame. The actual solutions 
$q_{j}(t)=e^{i\omega t}u_{j}(t)$ on the sphere can be obtained 
\textit{a posteriori} by the following inverse transformation:
\[
v_{j}(t)=\frac{1}{1+\left\vert q_{j}\right\vert ^{2}}\left(  2q_{j},\left\vert
q_{j}\right\vert ^{2}-1\right)  \in\mathbb{C\times R}\text{.}%
\]
To continue the numerical solutions in the plane, we use the augmented system
of equations
\begin{equation}
\dot{u}_{j}=T\left(  \lambda_{1}-i\right)  \omega u_{j}-T\left(  \lambda
_{2}-i\right)  {\frac{\left(  1+\left\vert u_{j}\right\vert ^{2}\right)  ^{2}%
}{4}}\sum_{k\neq j}\left(  \frac{u_{j}-u_{k}}{\left\vert u_{j}-u_{k}%
\right\vert ^{2}}-\frac{u_{j}}{1+\left\vert u_{j}\right\vert ^{2}}\right)
~\text{.}%
\end{equation}
The purely imaginary eigenvalues of the polygonal equilibrium can have
multiplicity~2, as was the case for vortices in the plane with or without
central vortex. The corresponding primary families of periodic orbits are then
determined via the perturbation scheme used in Section~\ref{sec:plane} and in
Section~\ref{sec:center}. When transformed back onto the sphere, the
stationary solutions of these equations in a rotating frame correspond to
polygonal equilibria at varying latitude, \textit{i.e.}, at varying value of
the spatial coordinate $Z$. The radius of such equilibria, after stereographic
projection onto the plane, equals $r=\cot(\theta/2)$, which introduces the
angle $\theta$ as a parameter of the problem. Note that $r\rightarrow\infty$
when $\theta\rightarrow0$, which on the sphere corresponds to the polygonal
equilibrium approaching the North pole. When $\theta\rightarrow{\frac{\pi}{2}%
}$ then $r \rightarrow1$ (the equator on the sphere), and when $\theta
\rightarrow\pi$ then $r\rightarrow0$ (the South pole). In particular, when
periodic solutions approach the North pole then $r\rightarrow\infty$, which
limits their numerical continuation. A simple way to continue the family
through the North pole is to change the chart of the stereographic projection
to the South pole, so that the North pole then corresponds to $r=0$.
\begin{figure}[ptb]
\begin{center}
\resizebox{12.5cm}{!}{
\includegraphics{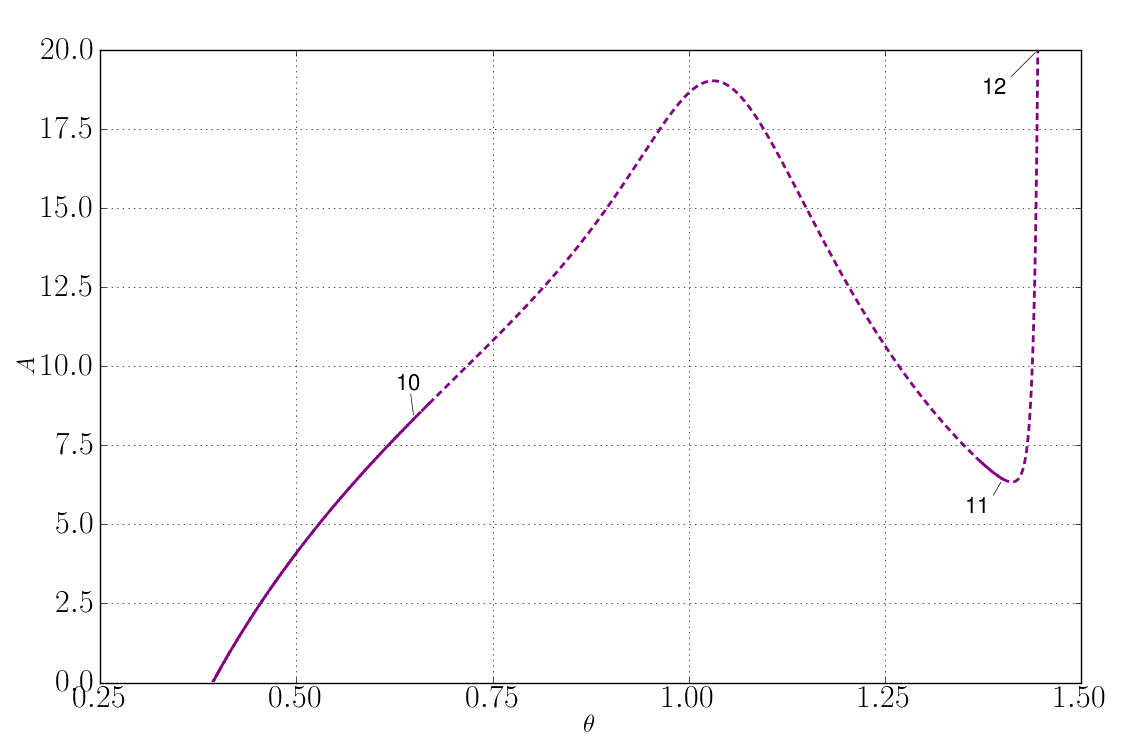}~}
\end{center}
\par
\vskip-.5cm\noindent\caption{ A continuum family of choreographies, for the
case of 5 vortices on a sphere. Solid curve sections represent stable or 
almost stable choreographies. The solutions labeled $10$, $11$, and $12$ are 
shown in Figures~\ref{fig08} and \ref{fig09}. Further data can be found in 
Table~1. }%
\label{fig07}%
\end{figure}

The presence of the problem parameter $\theta$ allows for additional
continuation schemes, as was the case in Section~\ref{sec:disk} for the disk,
and in Section~\ref{sec:center} for the case with a central vortex. For the
sphere we have found it useful to continue resonant orbits found along primary
families, keeping the resonance ratio fixed, as also done in
Section~\ref{sec:center}. In particular, such secondary continuation can be
used to determine a continuum of choreographies, as represented in
Figure~\ref{fig07} for the case of $5$ vortices. The particular family shown
there is obtained by first following one of the four families that bifurcate
from the relative equilibrium when $\theta=\pi/5$, namely one of the two
families that arise from the purely imaginary eigenvalue $2.78686i$ of
multiplicity~2. Along this primary family we can detect an unlimited number of
resonances that correspond to choreographies. Secondary continuation of one
such a resonant orbit, namely one with $T/T_{0}=-33:26$, keeping the resonance
fixed and allowing $\theta$ to vary, generates the family represented in
Figure~\ref{fig07}.

\begin{figure}[ptb]
\begin{center}
\resizebox{18.0cm}{!}{
\includegraphics{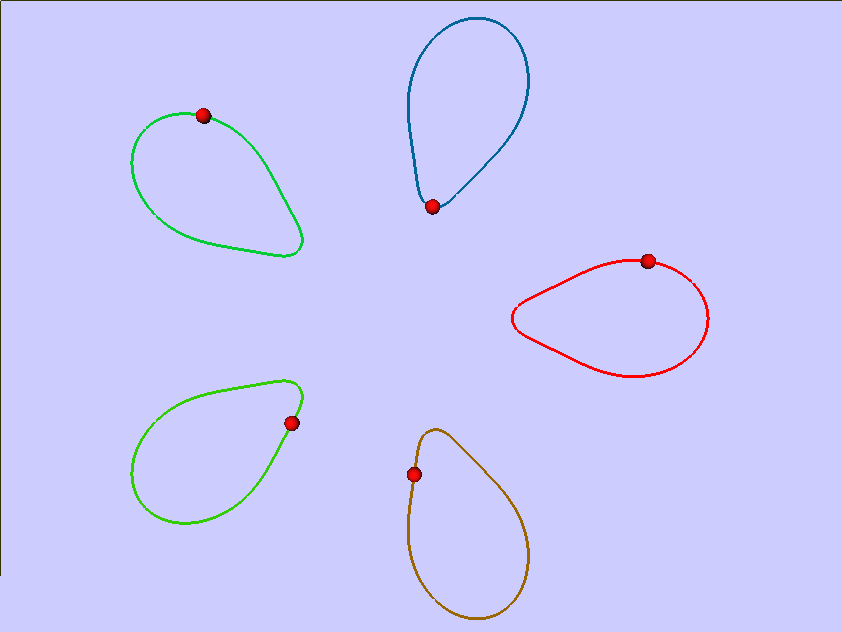}~
\includegraphics{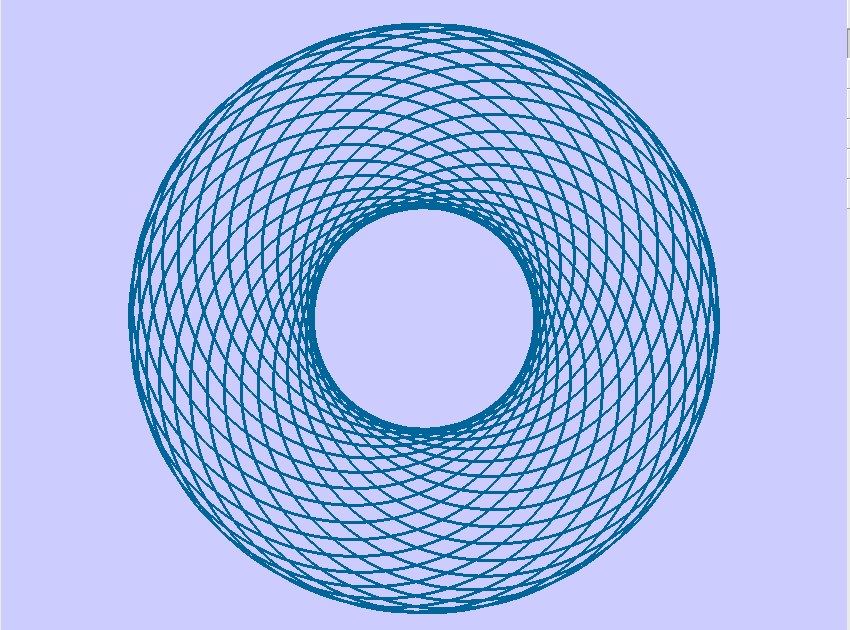} }
\end{center}
\par
\vskip-1.10cm\noindent
\par
\begin{center}
\resizebox{18.0cm}{!}{
\includegraphics{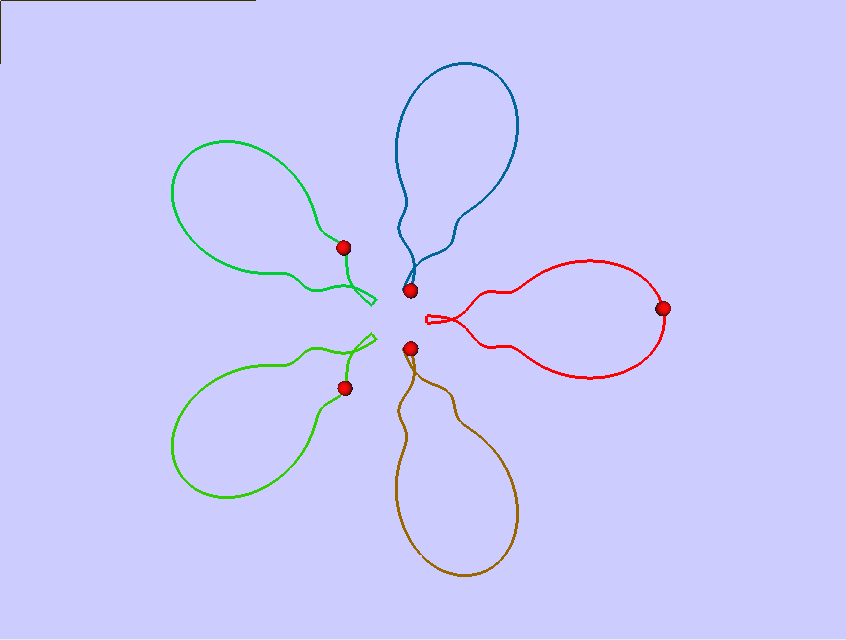}~
\includegraphics{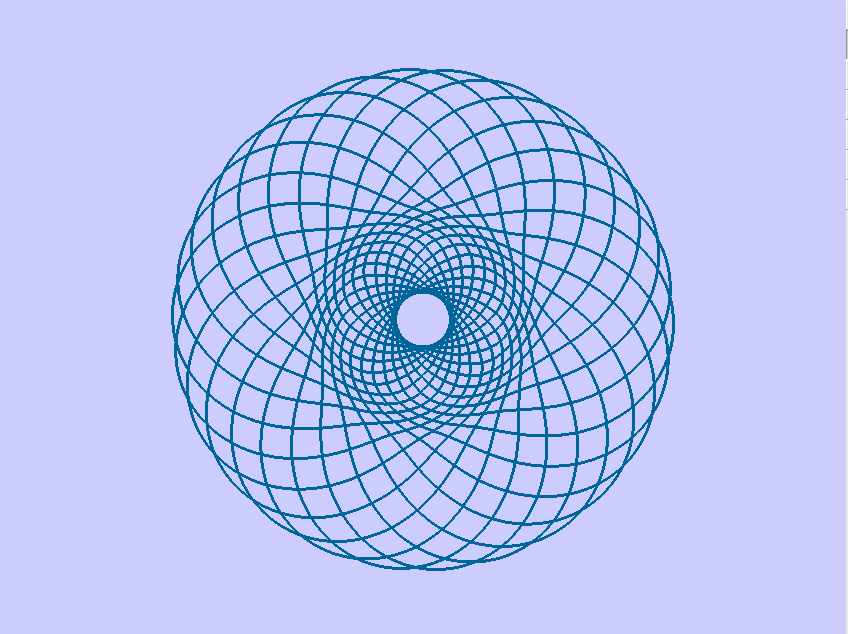} }
\end{center}
\par
\vskip-1.10cm\noindent
\par
\begin{center}
\resizebox{18.0cm}{!}{
\includegraphics{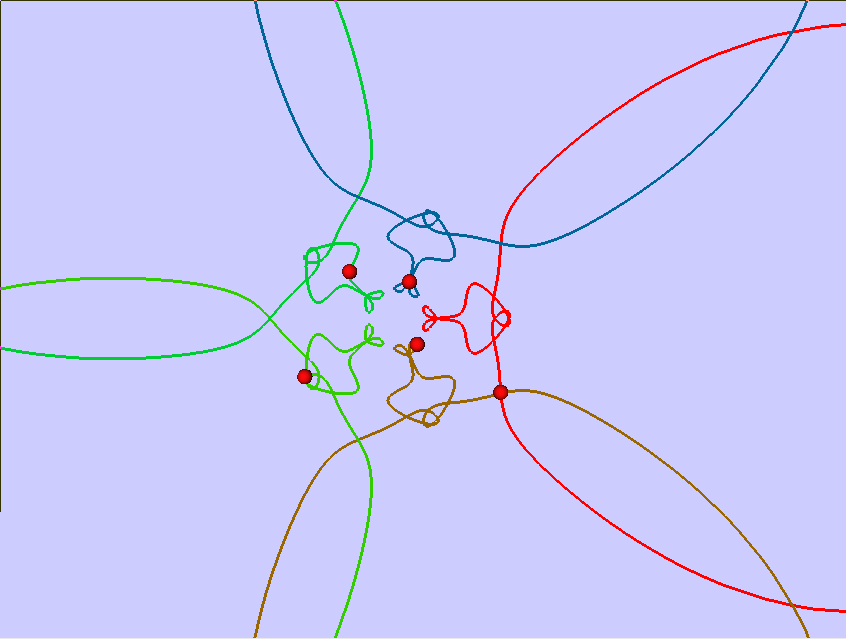}~
\includegraphics{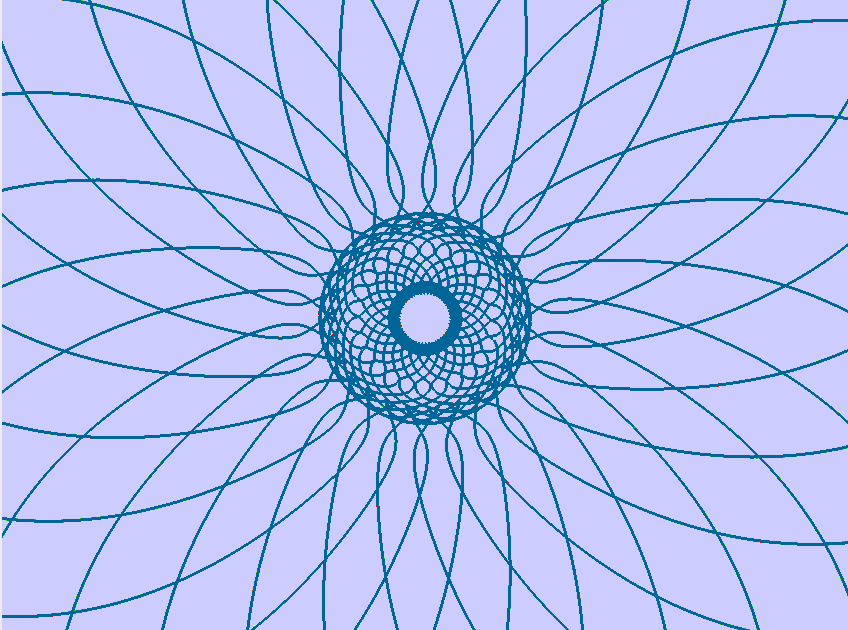} }
\end{center}
\par
\vskip-.5cm\noindent\caption{ Row-wise: Solutions $10$, $11$ and $12$ from the
continuum family of choreographies in Figure~\ref{fig07}, for 5 vortices on a
sphere, when stereographically projected onto the plane. The choreographies
are shown on the left in the rotating frame, and on the right in the inertial
frame. The orbits are drawn on a fixed scale, as a result of which the orbits
in the third row are shown only partially. Numerical data can be found in
Table~1. }%
\label{fig08}%
\end{figure}

\begin{figure}[ptb]
\begin{center}
\resizebox{18.0cm}{!}{
\includegraphics{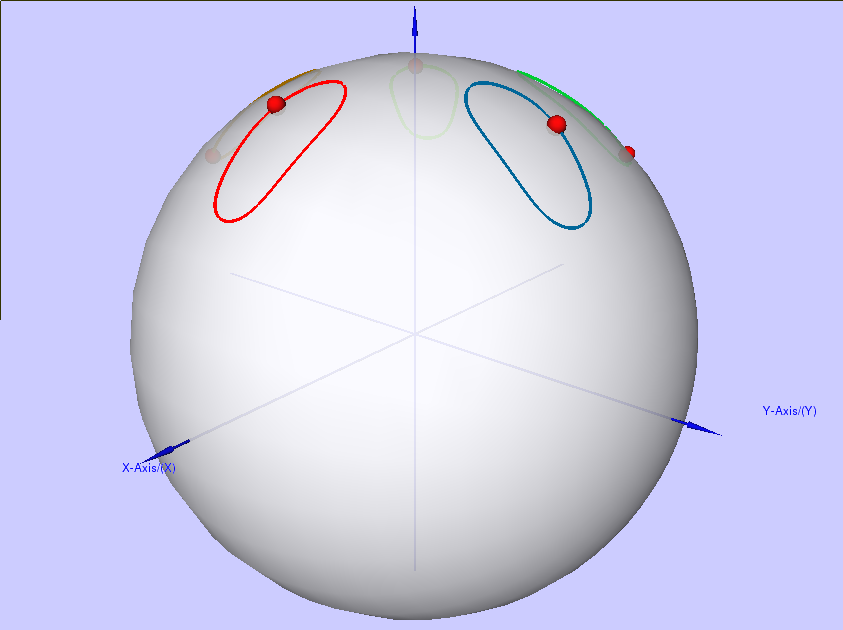}~
\includegraphics{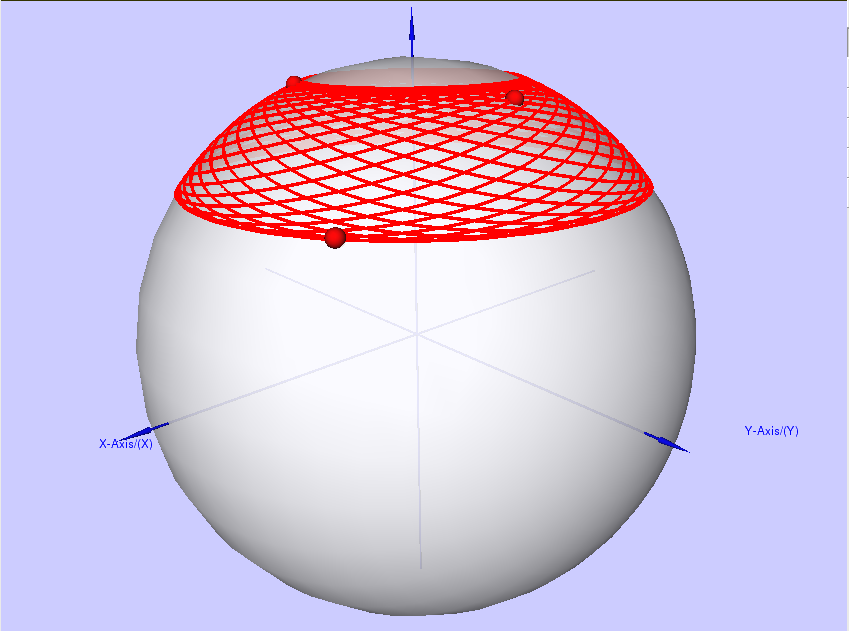} }
\end{center}
\par
\vskip-1.10cm\noindent
\par
\begin{center}
\resizebox{18.0cm}{!}{
\includegraphics{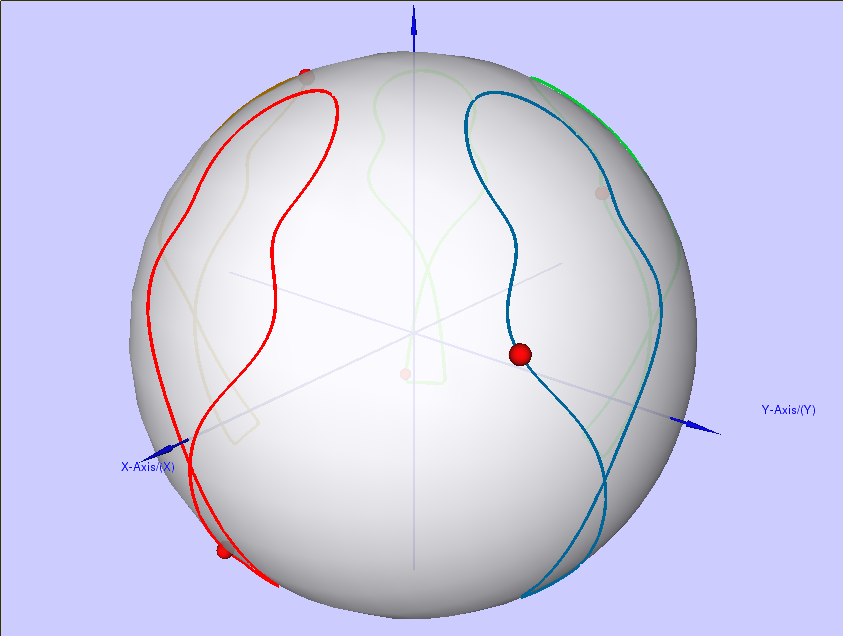}~
\includegraphics{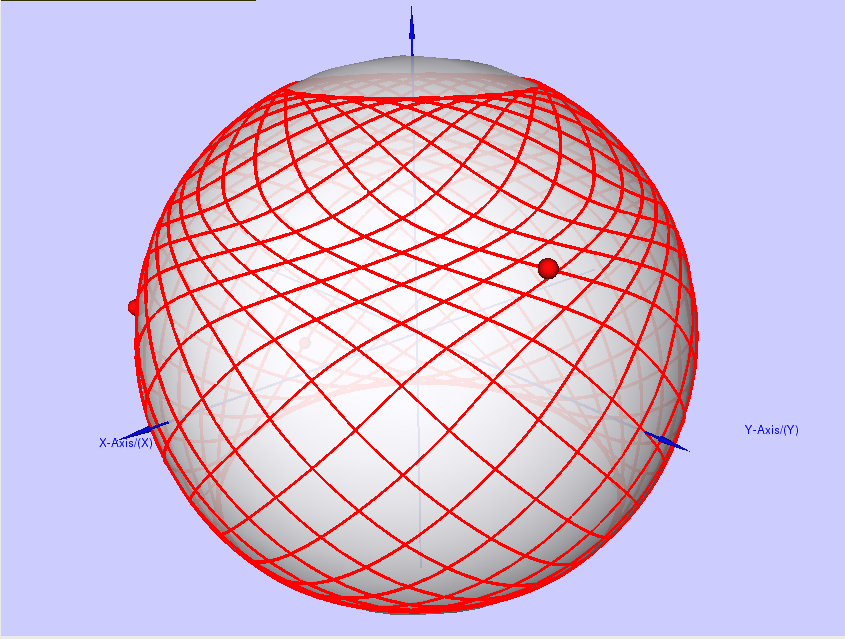} }
\end{center}
\par
\vskip-1.10cm\noindent
\par
\begin{center}
\resizebox{18.0cm}{!}{
\includegraphics{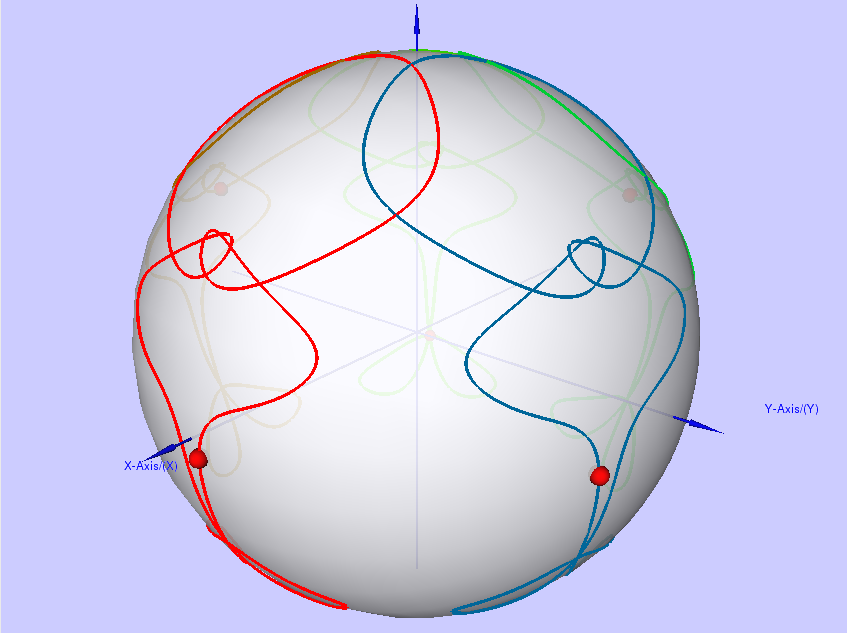}~
\includegraphics{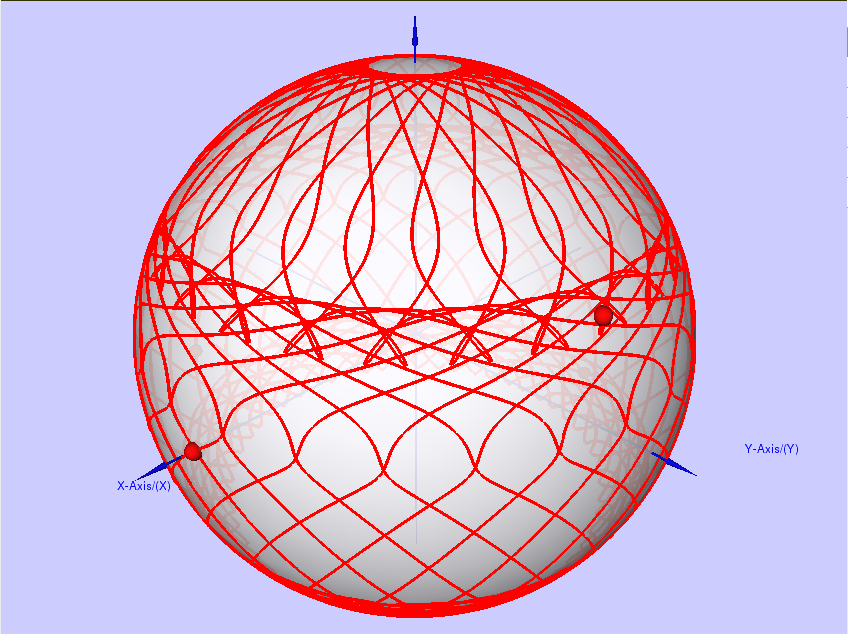} }
\end{center}
\par
\vskip-.5cm\noindent\caption{ Row-wise: Solutions $10$, $11$ and $12$ from the
continuum family of choreographies in Figure~\ref{fig07}, as shown in their
planar representation in Figure~\ref{fig08}. Here these orbits are shown after
being transformed back onto the sphere. The choreographies are shown on the
left in the rotating frame, and on the right in the inertial frame. Numerical
data can be found in Table~1. }%
\label{fig09}%
\end{figure}

The base period $T_{0}$ of the rotating frame is taken to be a function of
$\theta$, namely $T_{0}=2 \pi/\omega$, where $\omega=\omega(r)$, and
$r=\cot(\theta/2)$, with $\omega(r)$ defined in Equation~\ref{omega-sphere}.
We also note that the sign of $T_{0}$, and correspondingly the sign of a
resonance ratio, relates to the direction of rotation of the rotating frame.
The orbits labeled $10$, $11$, and $12$ in Figure~\ref{fig07}, of which the
first two are stable, are shown in Figure~\ref{fig08} (the planar
representation), and in Figure~\ref{fig09} (mapped onto the sphere); on the
left in the rotating frame, and on the right in the inertial frame. As follows
from the discussion above, all three have the same resonance ratio, namely,
$T/T_{0}=-33:26$. We note that the path to the family of choreographies in
Figure~\ref{fig07} could equally well have started from another polygonal
equilibrium for another value of $\theta$.

\begin{figure}[ptb]
\begin{center}
\resizebox{18.0cm}{!}{
\includegraphics{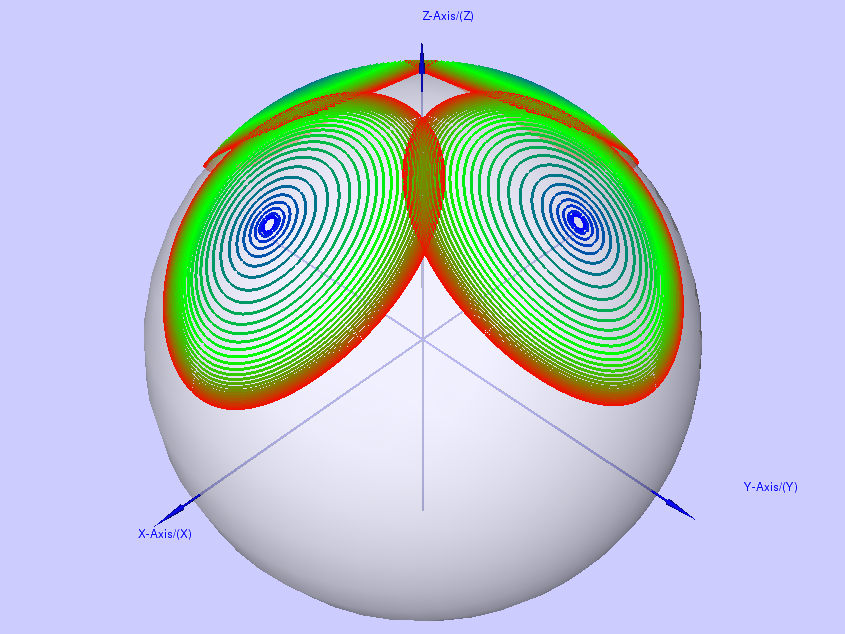}~
\includegraphics{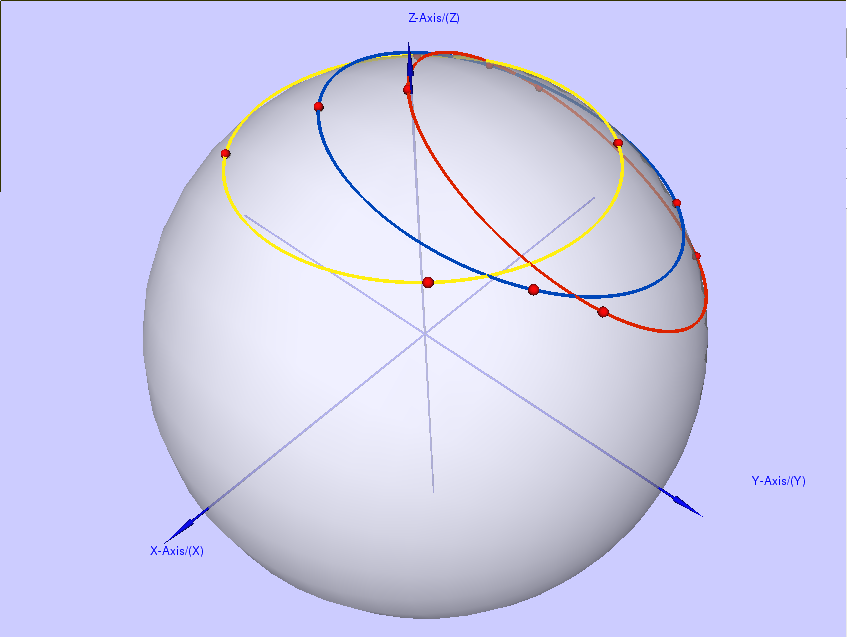} }
\end{center}
\par
\vskip-1.10cm\noindent
\par
\begin{center}
\resizebox{18.0cm}{!}{
\includegraphics{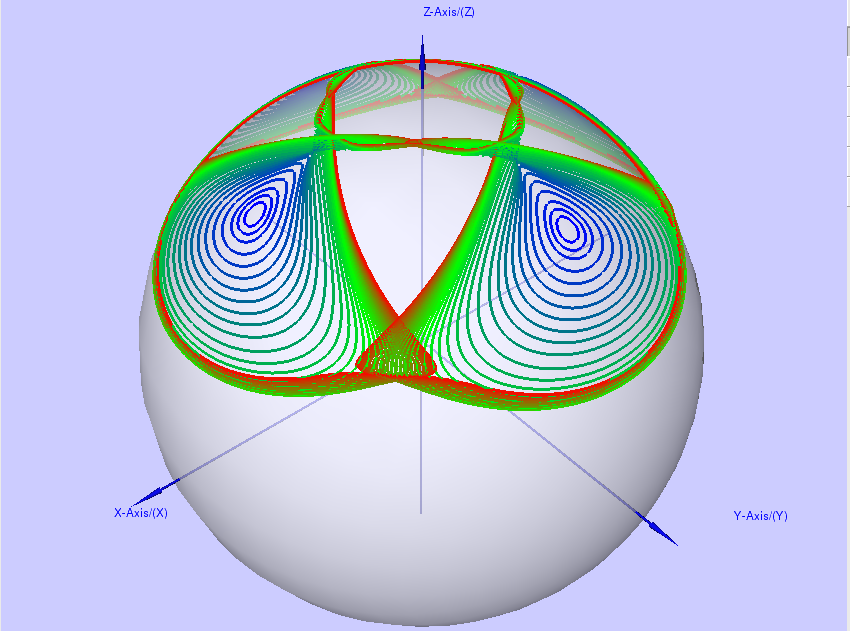}~
\includegraphics{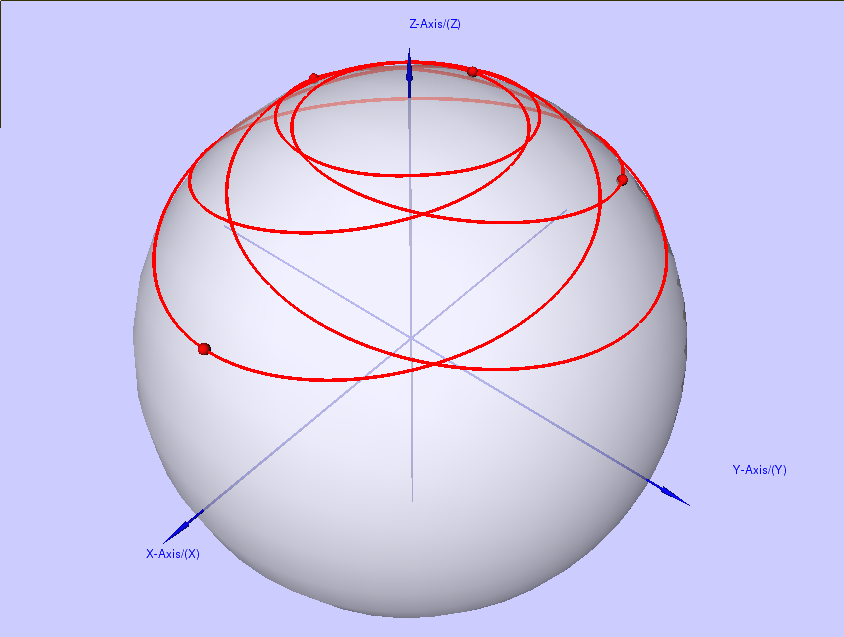} }
\end{center}
\par
\vskip-1.10cm\noindent
\par
\begin{center}
\resizebox{18.0cm}{!}{
\includegraphics{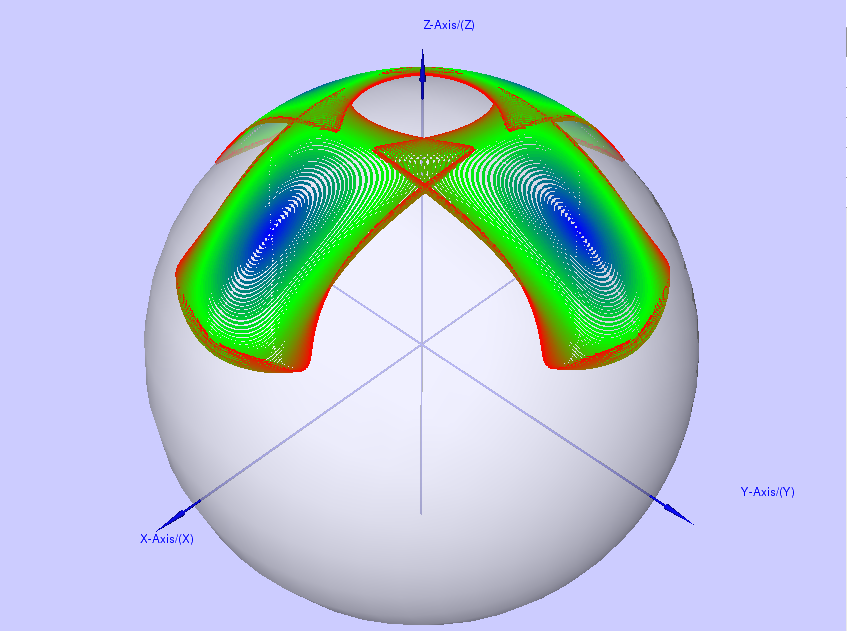}~
\includegraphics{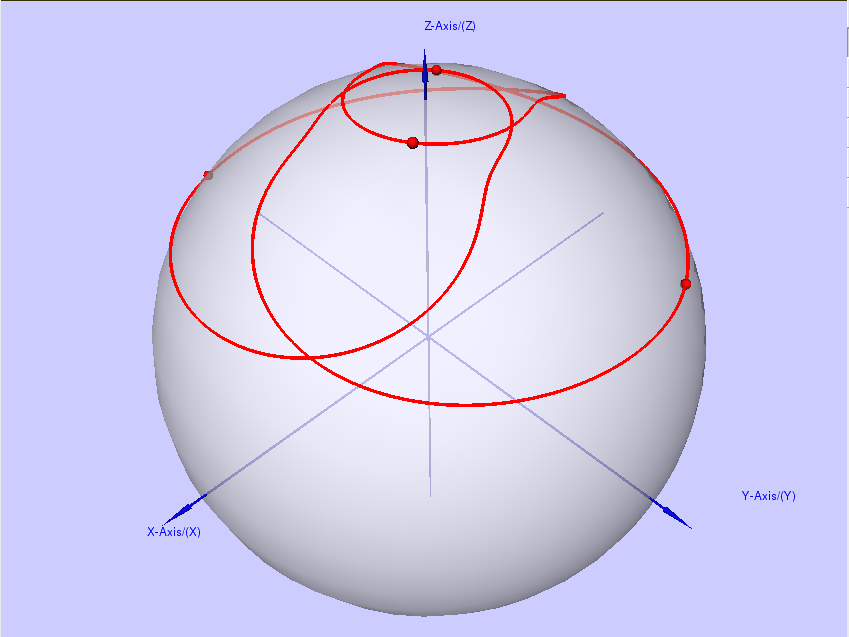} }
\end{center}
\par
\vskip-.5cm\noindent\caption{ Left columns: The three primary families that
bifurcate from the polygonal equilibrium for the case of $4$ vortices, with
$\theta=\pi/5$. Right columns: Selected choreographies that occur along the
primary families in the left column. }%
\label{fig10}%
\end{figure}

The three panels on the left in Figure~\ref{fig10} show a selection of orbits
along each of the three primary families that exist for the case of $n=4$
vortices. The families shown there arise from a polygonal equilibrium along
the northern hemisphere, namely for $\theta=\pi/5$, which has purely imaginary
eigenvalues $3.51246i$ (multiplicity $2$), and $2.84115i$ (multiplicity $1$),
as well as two zero eigenvalues. We note that a corresponding polygonal
equilibrium along the southern hemisphere, namely for $\theta=4\pi/5$, has the
same eigenvalues, but with period of the rotating frame of opposite sign,
namely $T_{0}=1.78883$ for the southern equilibrium, as opposed to
$T_{0}=-1.78883$ for the corresponding northern equilibrium.

The primary family shown in the top-left panel of Figure~\ref{fig10} is one of
the two families that arise from the double eigenvalue $3.51246i$. The orbits
along this family are degenerate, in the sense that the period $T$ remains
constant, with constant resonance ratio $T/T_{0}=-1$. In fact, all orbits of
this family correspond to choreographies, albeit trivial ones, as they
correspond to sliding a planar circular orbit along the sphere. Three such
circular choreographies are shown in the top-right column of
Figure~\ref{fig10}. The primary family shown in the center-left panel of
Figure~\ref{fig10} is the second family that arises from the double
eigenvalues $3.51246i$. Both families from this eigenvalue were determined via
the perturbation approach that was used earlier in this paper. One of the
countably infinite number of choreographies that exist along this second
family is shown in the center-right panel. The bottom-left panel of
Figure~\ref{fig10} shows orbits along the third family that arises from the
polygonal equilibrium, namely the family that arises from the multiplicity-1
eigenvalue $2.84115i$. This family can be computed directly, \textit{i.e.},
without use of the perturbation technique. Here also, one of the countably
infinite number of choreographies along this family is shown in the
bottom-right panel.

\begin{table}[]
\begin{center}
\begin{tabular}{|c|c|c|c||c|l|c|c|c|c|c|}
\hline
        &   &   &    &      &                &          &     &       \cr
Figure  &Row&Label&$n$&$T/T_0$&Parameter     &$T$       &$T_0$&Stable/\cr
        &   &   &    &      &                 &         &    &Unstable\cr
\hline
        &   &   &    &      &                 &         &         &   \cr
  2     &1  & 1 & 5  & 9:8  &                 & 3.53429 & 3.14159 & S \cr
        &   &   &    &      &                 &         &         &   \cr
(plane) &2  & 2 & 5  & 33:26& (none)          & 3.98741 & 3.14159 & S \cr
        &   &   &    &      &                 &         &         &   \cr
        &3  & 3 & 5  & 4:3  &                 & 4.18879 & 3.14159 & U \cr
        &   &   &    &      &                 &         &         &   \cr
\hline
        &   &   &    &      &                 &         &         &   \cr
  4     &1  & 4 & 5  & 13:6 & $R=1.99442$     & 6.78972 & 3.13372 & S \cr
        &   &   &    &      &                 &         &         &   \cr
(disk)  &2  & 5 & 5  & ~3:1 & $R=1.85172$     & 9.37522 & 3.12507 & U \cr
        &   &   &    &      &                 &         &         &   \cr
        &3  & 6 & 5  & 14:3 & $R=1.75556$     & 14.5297 & 3.11350 & S \cr
        &   &   &    &      &                 &         &         &   \cr
\hline
        &   &   &    &      &                 &         &         &   \cr
  6     &1  & 7 & 5+1& 11:4 & $\mu=1.85436$   & 4.48291 & 1.63015 & S \cr
        &   &   &    &      &                 &         &         &   \cr
 (with~~&2  & 8 & 5+1& 11:4 & $\mu=0.39014$   & 7.22918 & 2.62879 & U \cr
 center &   &   &    &      &                 &         &         &   \cr
 vortex)&3  & 9 & 5+1& 11:4 & $\mu=0.00500$   & 8.61784 & 3.13376 & U \cr
        &   &   &    &      &                 &         &         &   \cr
\hline
        &   &   &    &      &                 &         &         &   \cr
  8     &1  &10 & 5  &-33:26& $\theta=0.65$   & 1.83447 &-1.44534 & S \cr
        &   &   &    &      &                 &         &         &   \cr
(sphere)&2  &11 & 5  &-33:26& $\theta=1.40$   & 22.7821 &-17.9496 & S \cr
        &   &   &    &      &                 &         &         &   \cr
        &3  &12 & 5  &-33:26& $\theta=1.44466$& 31.1935 &-24.5767 & U \cr
        &   &   &    &      &                 &         &         &   \cr
\hline
\end{tabular}
\vskip0.cm\noindent
\caption{
Data for the orbits in Figures~1-8, where $n$ is the number of vortices, 
$T$ is the period of the orbit in the rotating frame,
and $T_0$ the period of the rotating frame.
The "parameter" column shows the value of a physical parameter, which
is absent in the case of vortices in the plane.  
Orbits are designated as ``Stable'' (S) when the largest Floquet multiplier 
is very close to the unit circle in the complex plane; otherwise they are 
listed as ``Unstable'' (U).
}
\end{center}
\end{table}

\section{Conclusions}


We have presented a systematic approach to the determination of choreographies
of the $n$-vortex problem in the plane, in a disk and on a sphere. Our
approach is based on the symmetries of the Lyapunov families of periodic
orbits that bifurcate from the polygonal relative equilibrium in a rotating
frame. The approach is similar to the one that we have used to find
choreographies for the $n$-body problem \cite{CaDoGa}. We use highly accurate
boundary value continuation techniques with adaptive meshes to compute the
Lyapunov orbits, and we presented a small, but representative selection of the
infinite number of choreographies that exists along the Lyapunov families. A
choreography of the $n$-polygon in the plane is invariant under scaling,
because the equations are homogeneous in that case. However, for a disk and a
sphere, the $n$-vortex problem is not homogeneous. One of the interesting
features of our results is therefore that in these cases we can determine a
continuum of choreographies with fixed resonance ratio, where the
choreographies are not related by scaling. Representative orbits from such
continua are included in our presentation.
\section*{ACKNOWLEDGMENTS}
This research was also supported by NSERC (Canada) Grant N00138. R.C. was partially supported by UNAM-PAPIIT Grant IA102818.


\end{document}